\documentclass[journal]{IEEEtran}
\usepackage{cite}

\ifCLASSINFOpdf
\else
\fi
\usepackage[cmex10]{amsmath}

\usepackage{cancel}

\usepackage{subfigure}
\usepackage{booktabs}     

\usepackage{stfloats}
\usepackage{graphicx}
\usepackage{color}
\usepackage{placeins}
\usepackage{float}
\usepackage[colorlinks=true,linkcolor=black,citecolor=black,urlcolor=black]{hyperref}
\usepackage{tabularx,colortbl,multirow}

\usepackage{enumerate}

\usepackage{bm}

\newcommand\bs[1]{\bm{\mathrm{#1}}}

\renewcommand\eqref[1]{(\ref{#1})}

\usepackage{tikz}
	\usetikzlibrary{shapes,decorations,decorations.pathreplacing}
	\usetikzlibrary{arrows,shapes,positioning,calc,fadings,patterns}
	\usetikzlibrary{decorations,decorations.pathreplacing,decorations.pathmorphing}
\usepackage{tikz-3dplot}

\usepackage{booktabs}     
\usepackage{threeparttable}
\usepackage{amssymb}
\usepackage{xcolor}
\usepackage{amsmath}
\usepackage{threeparttable}
\DeclareMathAlphabet{\mathpzc}{OT1}{pzc}{m}{it}

\renewcommand\epsilon{\varepsilon} 

\begin{document}

\title{Computation of  High-Order  Electromagnetic Field Derivatives with FDTD and the Complex-Step Derivative Approximation}

\author{Kae-An Liu, ~\IEEEmembership{Student~Member,~IEEE}, Hans-Dieter Lang, ~\IEEEmembership{Member,~IEEE}  and Costas~D.~Sarris,~\IEEEmembership{Senior~Member,~IEEE}
\thanks{The authors are  with the Edward S. Rogers Sr. Department of Electrical and
Computer Engineering, University of Toronto, Toronto, ON M5S 3G4, Canada
(e-mail: costas.sarris@utoronto.ca). This research has been supported by the Natural Sciences and Engineering 
Research Council of Canada (NSERC) through a Discovery Grant.} 
}

\markboth{IEEE TRANSACTIONS ON ANTENNAS AND PROPAGATIONS}%
{}
\maketitle

\begin{abstract}
This paper introduces a new approach for the computation of electromagnetic field derivatives, up to any order, with respect to the material and 
geometric parameters of a given geometry, in a single Finite-Difference Time-Domain (FDTD) simulation.  The proposed method is based on embedding the complex-step derivative (CSD) approximation into the standard FDTD update equations. Being finite-difference free,  CSD  provides accurate derivative approximations even for very  small perturbations of the design parameters, unlike finite-difference approximations that are prone to  subtractive cancellation errors. The availability of accurate approximations of field derivatives with respect to design parameters enables studies such as sensitivity analysis of multiple objective functions (as derivatives of those can be derived from field derivatives via the chain rule), uncertainty quantification, as well as  multi-parametric modeling  and optimization of electromagnetic structures. The theory, FDTD implementation and applications of this technique are presented. 
\end{abstract}


\IEEEpeerreviewmaketitle

\section{Introduction}
Quantifying the influence of material and fabrication tolerances  has been a topic of interest since the early days of computer-aided analysis and design of microwave circuits and systems \cite{bandlersevioraims70, bandlersevioraellett, 
tiberio-sensitivity}.  This topic is even more significant today, with the ever increasing complexity of electromagnetic structures. Nevertheless, research on computational electromagnetic solvers has primarily focused on  enhancing their speed and accuracy for well-defined problems, rather than problems defined with some degree of uncertainty and the evaluation of their relevant sensitivities. 

Variations in the geometry specifications and the materials of a structure can be considered in a sensitivity 
analysis \cite{saltelli-sensitivitybook, Martins-review},  by computing derivatives of an output function of interest (in the following referred to 
as an ``objective function"), with respect to the design parameters. A standard approach for this analysis is offered  by the 
finite-difference method. For  an input parameter $\xi$, with nominal value $\xi_0$,  
and an objective function $F(\xi)$, the centered finite difference (CFD): 
\begin{equation}
\dfrac{F(\xi_0 + h) - F(\xi_0 - h)}{2h}  = \dfrac{\partial F}{\partial \xi}(\xi_0)  + \dfrac{h^2}{6}\dfrac{\partial^3F}{
\partial \xi^3}(\xi_0) + \mathcal{O}(h^4)
\label{CFD}
\end{equation}
is a second-order accurate approximation to the derivative of $F$ with respect to $\xi$ (i.e. the leading error term 
of this approximation is proportional to $h^2$ and the notation $\mathcal{O}(h^4)$ implies that the rest of the terms are proportional to $h^4$ or higher powers of $h$).  This expression relies on iterative calls of a solver to calculate the objective function at points $\xi_0 \pm h$ of the parameter $\xi$, without any modification of the solver itself. However, as $h$ decreases, the absolute value of the difference  
$| F(\xi + h) - F(\xi - h) |$ may become smaller than machine precision (even if $F$ is computed analytically), or the error floor of the numerical method used to determine $F$ (if $F$ is numerically computed, which is the most practically interesting case). Then, the error of (\ref{CFD}) diverges and this approximation ceases to be useful.  This issue becomes even more important in the context 
of the Finite-Difference Time-Domain (FDTD) method, where small geometric and material perturbations are desired, to limit 
associated numerical dispersion errors.  Moreover, CFDs require two simulations per parameter, accumulating computational overhead when sensitivities with respect to  multiple parameters are considered.   

To address this latter issue rather than the former, adjoint variable methods (AVM), originally formulated for sensitivity analysis of electric circuits  \cite{director-avm},  have been implemented in frequency and time-domain numerical electromagnetic techniques \cite{nikolova-avm}. The key advantage of AVM is that it computes first-order sensitivities of an objective function with respect to multiple parameters, with just one additional simulation of the so-called adjoint problem. To do so  though, the adjoint solutions of $N$ perturbed problems with respect to each of the parameters are approximated by the adjoint solution to the unperturbed problem, hence 
further compromising the accuracy of the method. In terms of FDTD, \cite{nikolova-avm-FDTD} presented a  wave equation based formulation, which has been  recently extended to the computation of second-order sensitivities of an objective function \cite{Bakr-avm}. Using the wave equation rather than the first-order system of Maxwell's equations, which are discretized in conventional 
FDTD, is a drawback of these methods; a first step towards addressing this  was recently presented in \cite{bakr-avm-FDTD}.

To overcome the subtractive cancellation errors in finite-difference methods, first and higher-order derivative approximations 
can be computed by considering imaginary step perturbations $jh$ instead of the real steps employed in finite-difference expressions such as  (\ref{CFD}). Indeed, for  a real function $F$, the Taylor expansion is: 
\begin{multline}
F(\xi_0+jh)  = F(\xi_0) + jh \dfrac{\partial F}{\partial \xi}(\xi_0) - \dfrac{h^2}{2}
\dfrac{\partial^2 F}{\partial \xi^2}(\xi_0)  \\
-  j \dfrac{h^3}{6}  \dfrac{\partial^3 F}{\partial \xi^3}(\xi_0) + \mathcal{O}(h^4).
\label{complexstep-taylor}
\end{multline}
Evidently, (\ref{complexstep-taylor}) is complex-valued, whose real part leads to the approximation: 
\begin{equation}
\text{Re} \left\{F(\xi_0+jh) \right\} = F(\xi_0)   -  \dfrac{h^2}{2} \dfrac{\partial^2 F}{\partial \xi^2}(\xi_0) + 
 \mathcal{O}(h^4)
\label{realCSD}
\end{equation}
and its imaginary part leads to: 
\begin{equation}
\dfrac{\text{Im} \left\{ F(\xi_0+jh) \right\} }{h} =\dfrac{\partial F}{\partial \xi}(\xi_0)   - \dfrac{h^2}{6} \dfrac{\partial^3 F}{\partial \xi^3}(\xi_0)  +  \mathcal{O}(h^4)
\label{CSD}
\end{equation}
Therefore, the complex-step perturbation of a real function provides a second-order accurate approximation to the unperturbed value of the function as its real part and the derivative of the function with respect to the perturbed
parameter as its imaginary part.  Notably, (\ref{CSD}) is a finite-difference free (hence, free of the associated subtraction errors) approximation  called the complex-step derivative approximation (CSD \cite{martins-csd}).

It was recently pointed out \cite{sarris-lang-ims2015} that (\ref{CSD}) can be  directly embedded in numerical electromagnetic methods, such as FDTD, to enable the direct calculation of  electromagnetic field derivatives by 
a standard FDTD code.  The approach is simple: if the design parameters of interest are perturbed  by an imaginary step $j h $, the  fields $\bs{E}^n$, $\bs{H}^{n+1/2}$ produced by the update equations  at the $n$-th time step become complex; then $\text{Re} \left\{\bs{E}^n\right\}$, $\text{Re} \left\{\bs{H}^{n+1/2} \right\}$ are the fields of the unperturbed problem, and  $\text{Im} \left\{ \bs{E}^n \right\}/h$, $\text{Im} \left\{ 
\bs{H}^{n+1/2} \right\}/h$ are their derivatives with respect to the parameter that is perturbed by the imaginary step. 
Moreover, \cite{liu-sarris-ims2017} introduced an FDTD-based method to compute partial and high-order derivatives with respect to multiple parameters, using an  augmented version of CSD with multiple imaginary dimensions, the multi-complex step  derivative (MCSD) approximation \cite{Lantoine-mcsd}. 
\begin{figure*}[]
	\centering
		\subfigure[]{\includegraphics[width=8cm]{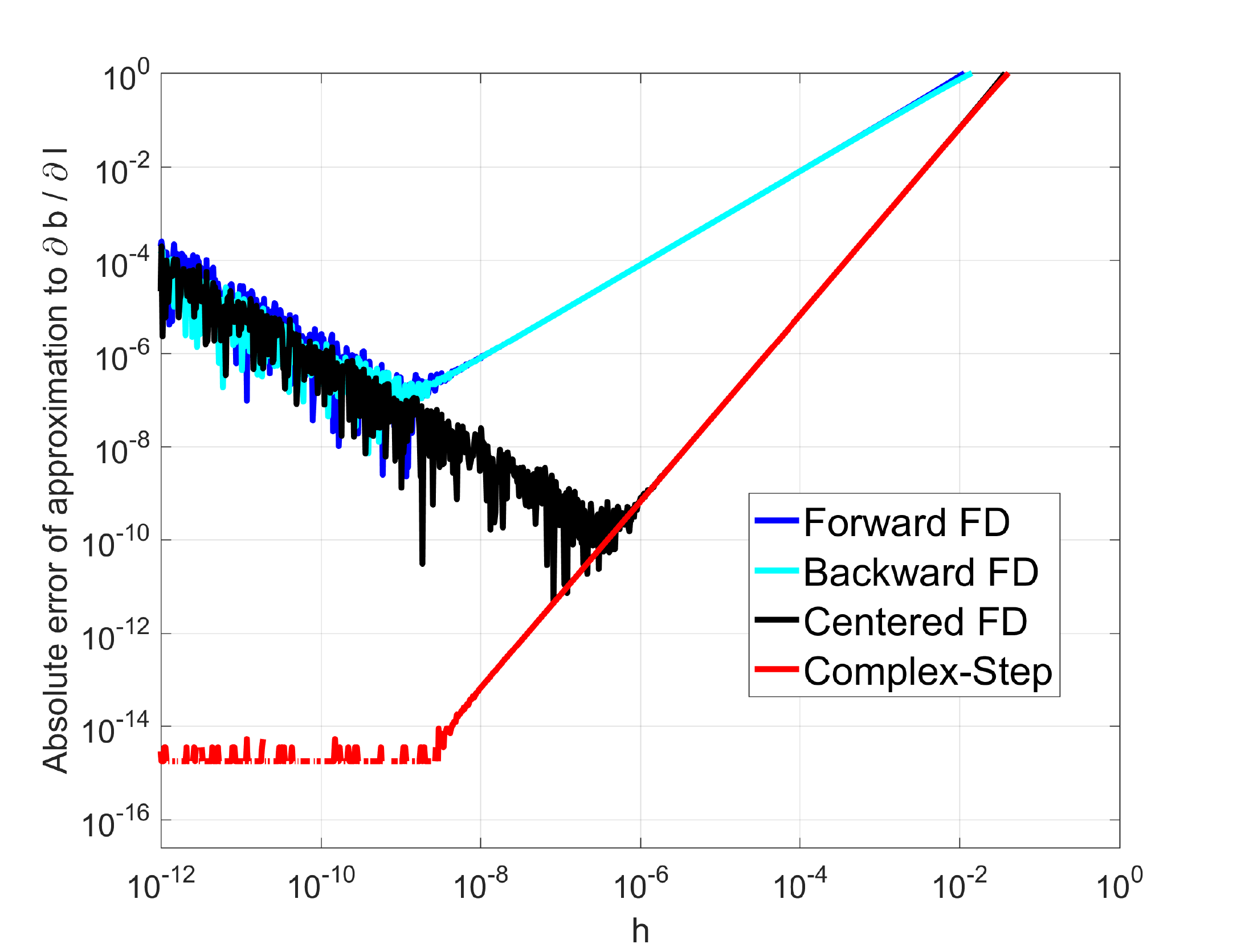}\label{fig:approxerrors}}
		\hspace{1cm}
		 \subfigure[]{\includegraphics[width=8cm]{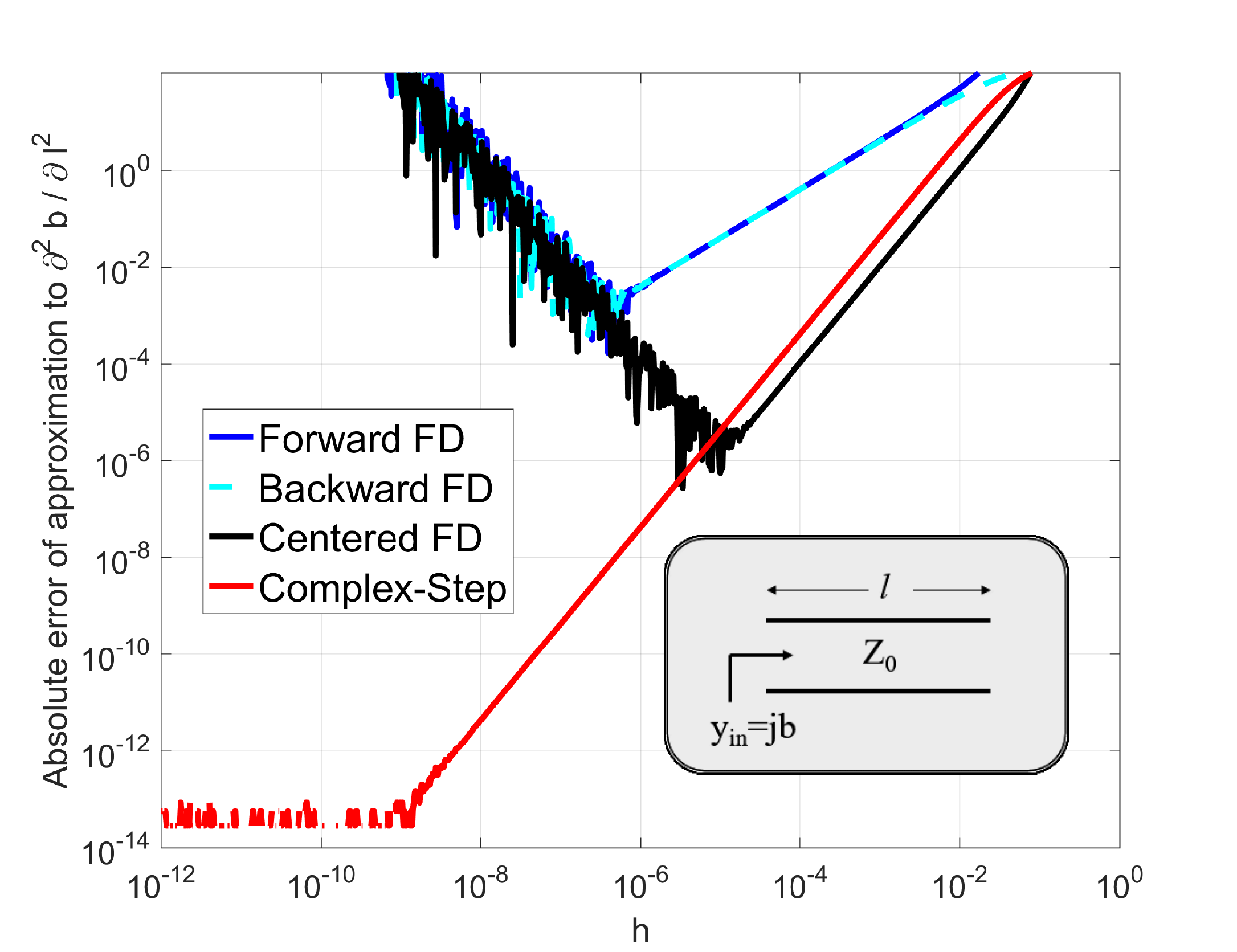} \label{fig:approxerrors2}}
		\caption{Absolute error in the approximation of the analytical value of (a) $\partial b / \partial l$ in 
(\ref{example_an}) and (b) $\partial^2 b / \partial l^2$ in (\ref{secondorder}) for forward, backward, and centered finite-differences as well as the complex-step derivative (CSD) approximation.}
		\label{filtersens}
\end{figure*}

There is ample motivation for further exploring the potential of this approach. First, embedding the CSD and MCSD 
in existing FDTD codes is straightforward. Second, such a combination allows for the computation of field derivatives on the fly; 
sensitivities of any field-based objective function can be found by post-processing, simply applying the chain 
rule. Third, a method that allows for the computation of any field derivative of any order is not just a tool for 
sensitivity analysis, but also a technique for parametric modeling \cite{hesthaven-parmodels,  liu-sarris-ursi2017} and uncertainty quantification \cite{hesthaven-TDUQ,  edwards-UQ-FDTD, austin-UQFDTD-mtt13, austin-UQchannel-tap13,  liu-sarris-iceaa2017}. Fourth, the high accuracy that complex step methods retain for small perturbations (unlike  finite-difference based methods) is particularly important for techniques such as FDTD, where geometric 
and material perturbations (typically modeled by stretching or squeezing Yee cells) can contribute to numerical dispersion errors. 

\textcolor{black}{
This paper builds on work reported in \cite{sarris-lang-ims2015, liu-sarris-ims2017, liu-sarris-ursi2017, liu-sarris-iceaa2017} to provide a comprehensive  presentation of general MCSD approximations implemented in FDTD, to formulate an accurate means of 
calculating field derivatives, {\em along} with the field solution to a problem. We extensively discuss the small, yet necessary  modifications of standard FDTD codes to embed MCSD approximations. The computational complexity and overhead of this method, with respect to 
the standard FDTD,  are thoroughly studied. We evaluate the accuracy of the method in a cavity case study with analytically known field derivatives and  we compute high-order field derivatives in a 3-D microwave circuit example. Simple as it is, this example illustrates how 
the subtractive cancellation errors of standard finite-difference methods limit their ability to compute second, third and higher-order derivatives. With high-order derivatives available, parameteric expressions of electromagnetic fields and field-based functions,
such as scattering parameters, are derived via Taylor expansions, demonstrating an important application of the proposed method 
to full-wave analysis based parametric modeling for electromagnetic design. 
}
\section{The Complex-Step Derivative Approximation}
A brief overview of the CSD and MCSD approximations is provided in this section, for the sake of completeness. 
We do so through specific examples.  Further details and notes on the historical origins of this approximation can be found in \cite{martins-csd, Lantoine-mcsd, Moler-numdiff}.
\subsection{Single complex-step derivative approximation}
Let us consider the normalized  input admittance of an open-ended, lossless transmission-line stub of 
length $l$: 
\begin{equation}
y_\mathrm{in} = j \tan   \dfrac{2 \pi l}{\lambda} 
\equiv  j b
 \label{stubyin}
\end{equation}
At $l/\lambda = 0.125$, the normalized susceptance $b=1$. The sensitivity of $b$ with respect  to the stub length can be expressed as: 
\begin{equation}
\dfrac{\partial b}{\partial l} = \dfrac{2 \pi}{\lambda} \dfrac{1}{\cos^2 \dfrac{2 \pi l }{\lambda} }
= \dfrac{2 \pi}{\lambda} \dfrac{1}{\cos^2 \dfrac{\pi }{4}} = \dfrac{4 \pi}{\lambda}
\label{example_an}
\end{equation}
The CSD approximation to this sensitivity is: 
\begin{equation}
\dfrac{\partial b}{\partial l} \approx   \text{Im} 
\dfrac{\tan  \dfrac{2 \pi (l+jh)}{\lambda} }{h} 
\label{example_csd}
\end{equation}
Forward, backward and centered difference approximations can also be used to find this sensitivity as 
follows: 
\begin{align}
\dfrac{\partial b}{\partial l} &\approx  
\dfrac{\tan  \dfrac{2 \pi (l+h)}{\lambda}  -\tan  \dfrac{2 \pi l}{\lambda}  }{h}
\hspace{.5cm}\text{(forward)}\\
\dfrac{\partial b}{\partial l} &\approx  
\dfrac{\tan   \dfrac{2 \pi l}{\lambda}  -\tan \dfrac{2 \pi (l-h)}{\lambda}  }{h}
\hspace{.5cm}\text{(backward)} \\
\dfrac{\partial b}{\partial l} &\approx  
\dfrac{\tan  \dfrac{2 \pi (l+h)}{\lambda}  -\tan   \dfrac{2 \pi (l-h)}{\lambda}  }{2h}
\hspace{.5cm}\text{(centered)} 
\end{align}
Normalizing $\lambda=1$, based on the analytical expression (\ref{example_an}), the errors of all four 
approximations can be found, as a function of $h$, for $l=0.125$. The results are shown in Fig. \ref{fig:approxerrors}. These confirm that 
CFD and CSD are second-order accurate with initially identical error performance, as also shown by the expressions (\ref{CFD}) and (\ref{CSD}). Yet,  the accuracy of CFD is clearly compromised by the subtraction error  of finite differences, while the accuracy of CSD is practically as good as that of the analytical solution up to machine precision. 

\subsection{Complex-step derivative approximation of second-order derivatives}
To derive a complex-step approximation to a second-order derivative, the bi-complex numbers are introduced. In 
general,  $n$-complex numbers $z \in \mathbb{C}^n $ are defined recursively as follows: 
\begin{equation}
\mathbb{C}^n := \{ z_1 + z_2 j_n | z_1, z_2 \in \mathbb{C}^{n-1} \} 
\label{ncomplexdef}
\end{equation}
For the bi-complex case ($n=2$), $j_2$ is the imaginary unit of a second imaginary dimension that is added to the one 
of the standard $\mathbb{C}$, for which the imaginary unit $j$ is written as $j_1$ in this generalized notation. 
For these two new imaginary units, and for any additional ones recursively introduced per (\ref{ncomplexdef}), 
\begin{subequations}
\begin{align}
j_1j_1&=j_2j_2=-1\\
j_1j_2&=j_2j_1\neq-1
\end{align}
\end{subequations}
Then, building a bi-complex number from the real and imaginary parts of the standard complex numbers leads to: 
\begin{subequations}
\begin{align}
z & =  z_1 + j_2 z_2 \\
  &  = \left( x_1 + j_1 y_1 \right) + j_2 \left( x_2 + j_1 y_2 \right)  \\
 & =   x_1 + j_1 y_1 + j_2 x_2 + j_1 j_2 y_2 
\end{align}
\end{subequations}
with: $\text{Re}(z) = x_1$, $\text{Im}_1 (z) = y_1$, $\text{Im}_2 (z) = x_2$, $\text{Im}_{12} (z) = y_2$. Note that 
there are three types of imaginary parts separately defined here, to distinguish $j_1$, $j_2$ and $j_1 j_2$
terms. Then,  a  bi-complex perturbation $\left( j_1 + j_2 \right)h$ can lead to a second-order derivative approximation 
as indicated by re-visiting the Taylor expansion (\ref{complexstep-taylor}): 
\begin{multline}
F \left( \xi_0 +  \left( j_1 + j_2  \right) h \right)  = F(\xi_0) + \left(  j_1 + j_2 \right) h 
\dfrac{\partial F}{\partial \xi}(\xi_0)  \\[.3cm] - 2 \left(1 - j_1 j_2 \right) \dfrac{h^2}{2}
\dfrac{\partial^2 F}{\partial \xi^2}(\xi_0)  
-  4 \left( j_1 + j_2 \right) \dfrac{h^3}{6}  \dfrac{\partial^3 F}{\partial \xi^3}(\xi_0) \\[.3cm]
 + 8 \left( 1 - j_1 j_2 \right) \dfrac{h^4}{24}  \dfrac{\partial^4 F}{\partial \xi^4}(\xi_0) 
 + \mathcal{O}\left( h^5 \right) 
\label{complexstep-taylor2}
\end{multline}
Letting $\text{Im}_{12}$ denote the term preceded by $j_1 j_2$, (\ref{complexstep-taylor2}) leads to the 
approximation (in addition to approximations (\ref{realCSD}), (\ref{CSD})): 
\begin{multline}
\dfrac{\partial^2 F}{\partial \xi^2}(\xi_0)  =  \dfrac{\text{Im}_{12} \left\{ F\left( \xi_0 +  \left( j_1 + j_2  \right) h \right) \right\} }{h^2} 
\\[.3cm] 
\:  + \dfrac{h^2}{3} \dfrac{\partial^4 F}{\partial \xi^4}(\xi_0) +   \mathcal{O}(h^4)
\label{CSD2}
\end{multline}
Note that the leading error term is proportional to $h^2$. However, as shown in Fig. 1(b), the constant in front of $h^2$ is 
four times larger than the one in the centered finite-difference approximation given below: 
\begin{multline}
\dfrac{\partial^2 F}{\partial \xi^2}(\xi_0) =  \dfrac{F(\xi_0+h) - 2 F(\xi_0) + F(\xi_0-h)}{h^2} \\[.3cm]
 -   \dfrac{h^2}{12} \dfrac{\partial^4 F}{\partial \xi^4}(\xi_0)  +  \mathcal{O}(h^4)
\label{CFD2}
\end{multline}
These formulas are applied to approximate the second order derivative of the normalized susceptance $b$ in (\ref{stubyin}), which is:
\begin{equation}
\dfrac{\partial ^2 b}{\partial l^2} = \dfrac{8 \pi^2}{\lambda^2} \dfrac{\sin  \dfrac{2 \pi l }{\lambda} }
{\cos^3 \dfrac{2 \pi l }{\lambda}} = \dfrac{16 \pi^2}{\lambda^2}
\label{secondorder}
\end{equation}
for $l=0.125$, $\lambda = 1$.  
The bi-complex step approximation to this second-order derivative is: 
\begin{equation}
\dfrac{\partial^2 b}{\partial l^2} \approx \text{Im}_{12} \dfrac{\tan \dfrac{2 \pi (l+(j_1 + j_2)h)}{\lambda}}{h^2} 
\label{CSD2-example}
\end{equation}
Expanding the tangent, one can readily re-cast (\ref{CSD2-example}) in a conventional $\mathbb{C}$ notation: 
\begin{equation}
\dfrac{\partial^2 b}{\partial l^2} \approx \text{Im} \dfrac{\zeta_2 \zeta_3 - \zeta_1 \zeta_4}{\zeta_3^2 + \zeta_4^2} 
\label{CSD2-exampleB}
\end{equation}
where $\text{Im} \equiv \text{Im}_1$, and: 
\begin{equation}
\begin{array}{lclcl}
\zeta_1 &  = &  \sin{\alpha} \cosh^2{\beta} & + & \dfrac{j}{2} \cos{\alpha} \sinh{2\beta} \\[.3cm]
\zeta_2 &  = &  \dfrac{1}{2}\cos{\alpha} \sinh{2 \beta} & - & j \sin{\alpha} \sinh^2{\beta} \\[.3cm]
\zeta_3 &  = &  \cos{\alpha} \cosh^2{\beta} & -  & \dfrac{j}{2} \sin{\alpha} \sinh{2\beta} \\[.3cm]
\zeta_4 &  = &  - \dfrac{1}{2}\sin{\alpha} \sinh{2 \beta} & - & j \cos{\alpha} \sinh^2{\beta} \\[.3cm]
\end{array}
\label{CSD2-exampleC}
\end{equation}
with $\alpha =  2 \pi l / \lambda$,  $\beta = 2 \pi h / \lambda$. This approximation is compared to standard forward, backward and centered finite-differences in  Fig. \ref{fig:approxerrors2}. Notably, subtraction errors dominate and, indeed, destroy the accuracy of finite-difference approximations, to a greater degree than 
in the case of first-order derivatives (Fig.  \ref{fig:approxerrors}), with the breakpoint moving to $h=10^{-5}$.  On the other hand, the bi-complex step approximation is able to deliver accuracy that is practically equivalent to that of the analytical expression of the second-order derivative, as $h \rightarrow 0$. In general, the accuracy 
advantage of the complex step approximations increases with the order of the derivative under consideration. With many applications, such as optimization studies 
and parametric modeling, requiring the computation of second and higher-order derivatives, this observation further motivates the study of complex-step  approximations. In the following, their generalized form  is presented. 
\subsection{Multi-complex step approximation for arbitrary order derivatives}
The last example has paved the way for the introduction of generalized complex-step derivative approximations for an arbitrary number of parameters and derivatives up to any order, introduced as the multi-complex step derivative (MCSD) approximation in \cite{Lantoine-mcsd}. 
Let us consider a function $F\left( \bs{\xi} \right)$, $\bs{\xi} = \left[ \xi_1, \xi_2, \cdots \xi_N  \right]^T$ and a partial derivative: 
\begin{equation}
\dfrac{\partial^{m_1+m_2+\cdots m_N }F}{\partial \xi_1^{m_1} \partial \xi_2^{m_2} \cdots
\partial \xi_N^{m_N}}
\label{hoderiv}
\end{equation}
The MCSD approach is to use $m_1 + m_2 + \cdots m_N \equiv P$ imaginary dimensions, i.e. perform 
computations in $\mathbb{C}^{P}$, introducing $m_k$ imaginary perturbations of each parameter $\xi_k$, 
$k=1, \cdots N$, with magnitude $h_k$. Then, the MCSD approximation of (\ref{hoderiv}) becomes: 
\begin{equation}
\begin{array}{l}
\dfrac{\partial^{m_1+m_2+\cdots m_N}F}{\partial \xi_1^{m_1} \partial \xi_2^{m_2} \cdots
\partial \xi_N^{m_N}} \approx \\[.3cm]
\dfrac{\text{Im}_{m_1,\cdots, m_N} F\left( \xi_1 + h_1 \cdot  \hspace{-.1cm} \displaystyle \sum_{k=1}^{m_1} j_k , \cdots 
,\xi_N + h_N \cdot  \hspace{-.6cm}  \displaystyle \sum_{k=P-m_N+1}^P j_k \right)}{h_1^{m_1} h_2^{m_2} \cdots h_N^{m_N}}
\label{MCSD}
\end{array}
\end{equation}
In this expression, $\text{Im}_{m_1,\cdots, m_N}$ means that once all perturbations have been applied to 
$F$, the approximation is actually derived from the $j_{m_1} j_{m_2} \cdots j_{m_N}$ term, just as the 
second-order derivative approximation (\ref{CSD2-example}) was derived from the $j_1 j_2$ term. 

The general approximation (\ref{MCSD}) and its lower-dimension counterparts, presented through the examples 
of the previous subsections, can be applied to several numerical methods for the analysis, design and optimization 
of microwave circuits. In the following, we focus on how this approximation can be embedded into the 
FDTD technique, to compute arbitrary order field derivatives with respect to design parameters along with the 
full-wave solution to a given problem. 
\section{FDTD-based Computation of Field Derivatives using the Complex-Step Derivative Approximation}
\subsection{Formulation}
\begin{figure}[ht]
	\centering
		\includegraphics[width=8cm]{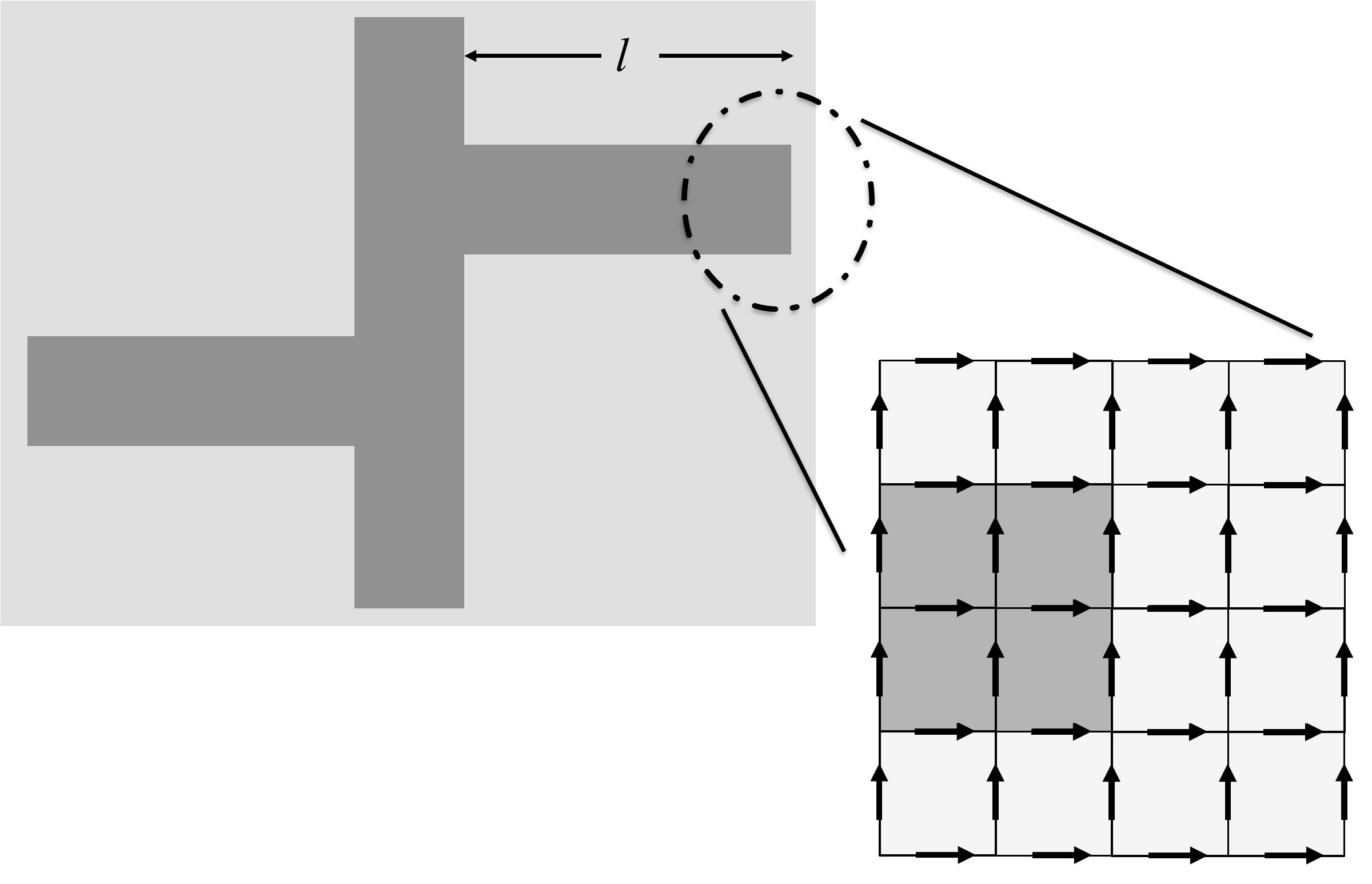}  
	\caption{Example of a microstrip circuit: field derivatives with respect to the stub length $l$ and the dielectric 
	permittivity, $\epsilon_r$, of the substrate are sought for. Dark grey areas indicate the microstrip segments; light grey areas correspond to 
	the substrate.}
	\label{msexample}
\end{figure}

The general strategy for the implementation of CSD in FDTD is to introduce complex perturbations into the standard
FDTD update equations. We explain this through an example of a microstrip segment (Fig. \ref{msexample}), 
printed on a dielectric substrate. Here, first and second order field derivatives with respect to the length $l$ of the segment and the relative dielectric permittivity  $\epsilon_r$ of the substrate are sought for.  

\subsubsection{Field derivatives with respect to material parameters}
To compute derivatives with respect to  $\epsilon_r$, two imaginary dimensions $j_1$ and $j_2$
are introduced and $\epsilon_r$ is set to  $\epsilon_r  + \left( j_1 + j_2\right)h_1$ in the FDTD 
update equations. For example, the update of the $x$-component of the electric field becomes: 
\begin{multline}
E_{i',m,k}^{x,n+1} =  E_{i',m,k}^{x,n}  +  \dfrac{\Delta t}{\epsilon_0 
\left( \epsilon_r  +  \left( j_1 + j_2 \right)  h_1 \right) } \times
 \\[.4cm]
 \Bigg( \dfrac{H_{i',m',k}^{z,n'}- H_{i',m'-1,k}^{z,n'}}{\Delta y} 
 -  \dfrac{H_{i',m,k'}^{y,n'} -
H_{i',m,k'-1}^{y,n'}}{\Delta z} \Bigg)
\label{mcsdupdate}
\end{multline} 
In (\ref{mcsdupdate}) and the FDTD update equations that follow $\left(i,m,k\right)$ are Yee cell indices, $n$ is a time-step index and  $p'=p+1/2$ ($p=i,m,k,n$) is the half cell/time-step offset  of space/time nodes in the FDTD mesh. Once a class of $n$-complex numbers is available (as discussed in the next subsection), this update equation is 
implemented in $\mathbb{C}^{2}$ and  the MCSD approximation
to the derivative $\partial^2 E_{i',m,k}^{x,n} / \partial \epsilon_r^2$ is simply: 
\begin{equation}
\dfrac{\partial^2 E_{i',m,k}^{x,n}}{\partial \epsilon_r^2} \approx 
\dfrac{\text{Im}_{12} E_{i',m,k}^{x,n}}{{h_1}^2}
\end{equation}
To gain further insights into this equation and how field derivatives are deduced from it, one can  decompose it into real and (multiple) imaginary parts. To that end, it is recognized that: 
\begin{equation}
\begin{array}{l}
\dfrac{\Delta t}{\epsilon_0 
\left( \epsilon_r  +  \left( j_1 + j_2 \right)  h_1 \right) } =  \\[.3cm]
\dfrac{\Delta t}{\epsilon_0 \epsilon_r \left(1 + 4 h_1^2\right)}
\dfrac{\left( \epsilon_r^2 + 2 h_1^2 \right) - ( j_1 + j_2 )h_1 \epsilon_r - j_1j_2 h_1^2 }
{\epsilon_r^2 + 4 h_1^2} \equiv  \\[.4cm]
\alpha^R + j_1 \alpha^{I_1} + j_2 \alpha^{I_2} + j_1 j_2 \alpha^{I_{12}}
\end{array}
\end{equation}
where the super-scripts $R$, $I_1$, $I_2$, $I_{12}$ indicate the real and multi-imaginary ($j_1$, $j_2$, 
$j_1j_2$) terms, respectively.  Then, the update equation for $\text{Im}_{12} E_{i',m,k}^{x,n}$
$\equiv$ $E_{i',m,k}^{I_{12},x,n}$ becomes: 
\begin{multline}
E_{i',m,k}^{I_{12}, x,n+1} =  E_{i',m,k}^{I_{12}, x,n}  \\[.5cm]
+ \alpha^{I_{12}} \Bigg( \dfrac{H_{i',m',k}^{R,z,n'} -
H_{i',m'-1,k}^{R,z,n'}}{\Delta y} -  \dfrac{H_{i',m,k'}^{R, y,n'}-
H_{i',m,k'-1}^{R, y,n'}}{\Delta z} \Bigg)   \\[.5cm]
+ \alpha^{I_2}
\Bigg(
\dfrac{H_{i',m',k}^{I_1,z,n'}- H_{i',m'-1,k}^{I_1,z,n'}}{\Delta y}   
 -  \dfrac{H_{i',m,k'}^{I_1, y,n'}-
H_{i',m,k'-1}^{I_1, y,n'}}{\Delta z} \Bigg) 
  \\[.5cm]
+  \alpha^{I_1}
 \Bigg(
\dfrac{H_{i',m',k}^{I_2,z,n'} -
H_{i',m'-1,k}^{I_2,z,n'}}{\Delta y}  
 -  \dfrac{H_{i',m,k'}^{I_2, y,n'}-
H_{i',m,k'-1}^{I_2, y,n'}}{\Delta z} \Bigg)   \\[.5cm]
 + \alpha^{R} \Bigg(
\dfrac{H_{i',m',k}^{I_{12},z,n'}- H_{i',m'-1,k}^{I_{12},z,n'}}{\Delta y}    
 -  \dfrac{H_{i',m,k'}^{I_{12}, y,n'}-
H_{i',m,k'-1}^{I_{12}, y,n'}}{\Delta z} \Bigg)
\label{update12}
\end{multline} 
Inspection of  (\ref{update12}) reveals that (\ref{mcsdupdate}) is equivalent to standard FDTD-type finite 
difference equations, which couple the real and imaginary fields. Note also that $|\alpha^R|$, $|\alpha^{I_1}|$, $|\alpha^{I_2}|$, $|\alpha^{I_{12}}|$ $\leq$ $\Delta t / \epsilon_0 \epsilon_r$.  So, if these finite differences 
appeared in standard FDTD, they would correspond to media with $\tilde{\epsilon}_r > \epsilon_r$. As a result,  these updates are subject to the standard Courant-Friedrichs-Lewy (CFL) stability condition of FDTD. Moreover, the memory requirements of MCSD-FDTD are greater than those of standard FDTD 
by a factor equal to $2^n$, for computations in $\mathbb{C}^n$. For $n=2$, for example, the real part of the fields is accompanied by the 
$I_1$, $I_2$ and $I_{12}$ terms. 
\subsubsection{Field derivatives with respect to geometric properties}
\label{geom-MCSD}
\begin{figure}[]
	\centering
		\includegraphics[width=8cm]{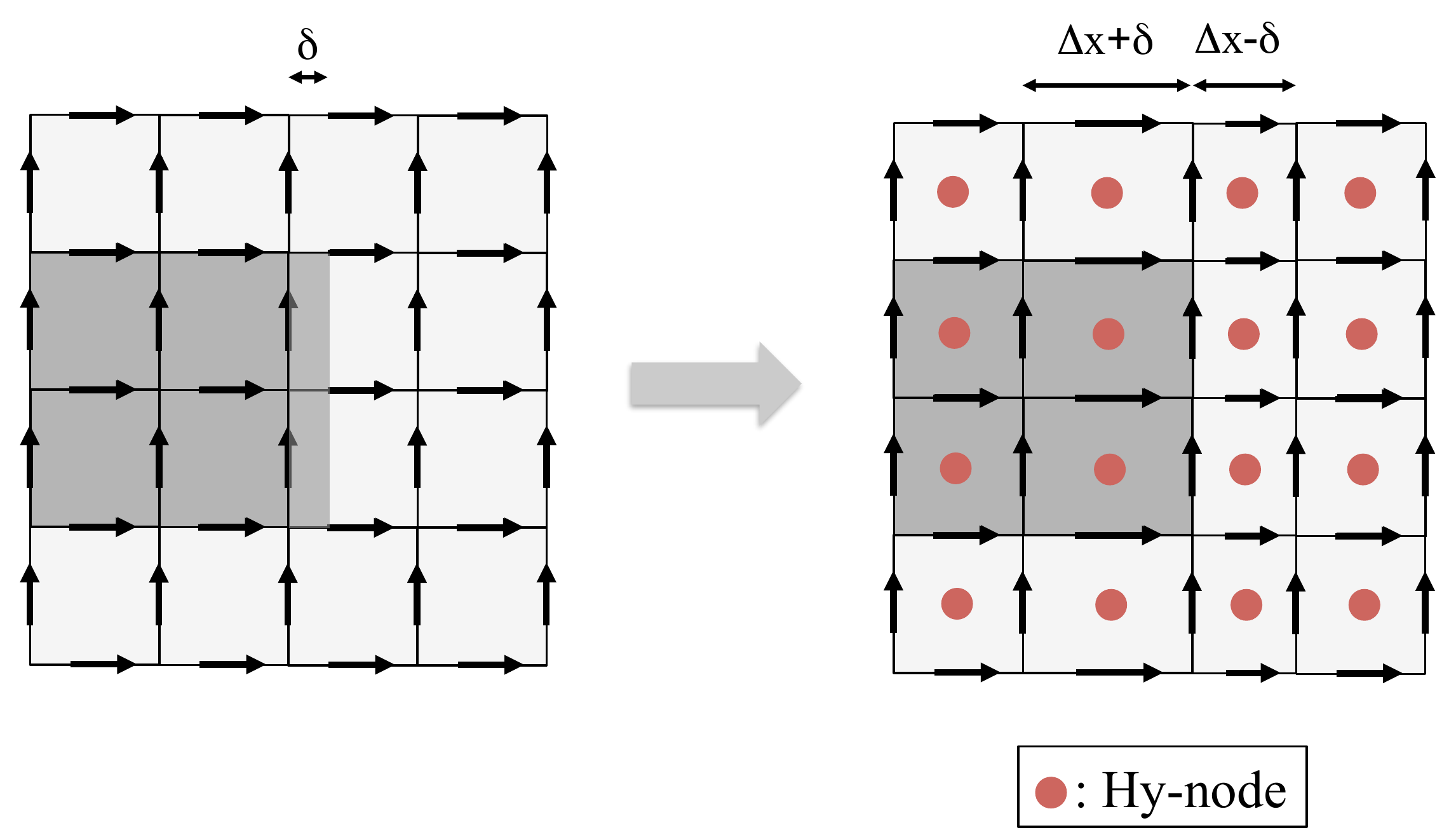}  
	\caption{FDTD mesh for the computation of field derivatives with respect to the length $l$ of the microstrip  stub
	in Fig. \ref{msexample}.  }
	\label{mcsdmesh}
\end{figure}

To derive the MCSD expressions for field derivatives with respect to the stub length $l$ in Fig. \ref{msexample},
the corresponding analysis with the finite-difference method is invoked.  In standard finite-difference sensitivity analysis, a common approach is to locally stretch or compress  the Yee cells by a geometric perturbation  $\delta$ of the length of the cells at the edge of the stub \cite{nikolova-avm-FDTD}, 
as shown in Fig. \ref{mcsdmesh}. 
Accordingly, the $H_y$ node shown in the figure can be updated as: 
\begin{multline}
H_{i',m,k'}^{y,n+1/2} =    H_{i',m,k'}^{y,n-1/2}  
+  \dfrac{\Delta t}{\mu_0 \Delta z} \left( E_{i',m,k+1}^{x,n}-  E_{i',m,k}^{x,n} \right) \\[.4cm]
-    \dfrac{\Delta t}{\mu_0 \Delta x^{-}} 
\left( E_{i+1,m,k'}^{z,n} - \cancelto{0} {E_{i,m,k'}^{z,n}} \right)
\end{multline}
where $\Delta x^{-} = \Delta x - \delta$. Likewise, update equations on stretched cells are derived by replacing 
$\Delta x$ by $\Delta x^{+} = \Delta x + \delta$. Translating this approach into the context of MCSD, the same 
form of equations is used, with  $\Delta x^{+,-} =   \Delta x \pm\left(  j_3  + j_4 \right) h_2$, where imaginary dimensions have been introduced through $j_{3,4}$, in addition to those corresponding to 
the derivatives with respect to $\epsilon_r$. In other words, there is a one-to-one mapping between the perturbations applied to a finite-difference based sensitivity analysis and MCSD.  

The relations between the derivatives of the electric field ${\textbf{E}}_{i,m,k}^n$   with respect to $l$, $\epsilon_r$ and the solution ${\tilde{\textbf{E}}}_{i,m,k}^n$ of the problem with the complex-step perturbations, 
are summarized here:
\begin{equation}
\begin{array}{l}
\text{Re} \left\{ \tilde{\textbf{E}}_{i,m,k}^n \: \right\}  \approx   {\textbf{E}}_{i,m,k}^n  \\[.2cm]
\text{Im}_{1} \left\{ \tilde{\textbf{E}}_{i,m,k}^n \: \right\}/h_1  \approx  \text{Im}_{2} \left\{ \tilde{\textbf{E}}_{i,m,k}^n \right\}/h_1  \approx  \partial \textbf{E}_{i,m,k}^n / \partial \epsilon_r \\[.2cm]
\text{Im}_{3} \left\{ \tilde{\textbf{E}}_{i,m,k}^n \: \right\}/h_2 \approx 
\text{Im}_{4} \left\{ \tilde{\textbf{E}}_{i,m,k}^n \: \right\}/h_2  \approx 
\partial \textbf{E}_{i,m,k}^n / \partial l  \\[.2cm]
\text{Im}_{12} \left\{ \tilde{\textbf{E}}_{i,m,k}^n \: \right\}/h_1^2 \approx
\partial^2 \textbf{E}_{i,m,k}^n / \partial \epsilon_r^2 
\\[.2cm] 
\text{Im}_{13} \left\{\tilde{\textbf{E}}_{i,m,k}^n \: \right\}/(h_1 h_2) \approx \partial^2 \textbf{E}_{i,m,k}^n / 
\partial \epsilon_r \partial l  \\[.2cm] 
\text{Im}_{34} \left\{ \tilde{\textbf{E}}_{i,m,k}^n \: \right\}/h_2^2 \approx \partial^2 \textbf{E}_{i,m,k}^n / \partial l^2   \\[.2cm]
\text{Im}_{123} \left\{ \tilde{\textbf{E}}_{i,m,k}^n \: \right\}/(h_1^2 h_2) \approx \partial^3 \textbf{E}_{i,m,k}^n / 
 \partial \epsilon_r^2 \partial l \\[.2cm]
\text{Im}_{134} \left\{ \tilde{\textbf{E}}_{i,m,k}^n \: \right\}/(h_1 h_2^2) \approx \partial^3 \textbf{E}_{i,m,k}^n / 
 \partial \epsilon_r \partial l ^2  \\[.2cm]
\text{Im}_{1234} \left\{ \tilde{\textbf{E}}_{i,m,k}^n \: \right\}/(h_1^2h_2^2) \approx \partial^4 \textbf{E}_{i,m,k}^n / \partial \epsilon_r^2 \partial l^2  
\end{array}
\label{mcsdoutput}
\end{equation}
Note that the real part of the field arrays contains the solution to the problem and the imaginary dimensions
contain derivatives from first to fourth order, including all entries to the Hessian matrix and all derivatives needed 
for a second-order Taylor series expansion of the fields with respect to ($l$, $\epsilon_r)$.  Hence, MCSD-FDTD is  more than a sensitivity analysis technique; it is a general framework for field derivative calculation, which can be 
readily employed in parametric modeling and uncertainty quantification studies. 
\subsection{FDTD Computations in $\mathbb{C}^n$}
A key step for embedding the MCSD approximation in existing FDTD codes is the implementation of a programming class for multi-complex numbers. With this class at hand, the FDTD arrays are defined as $n$-complex (with  $n$ chosen according to the maximum order of the derivatives that are needed) and perturbations of the imaginary part of any design parameter can be introduced, without any other change in the structure of the code. 

Since no standard multi-complex number data type is currently available,  a custom class, equipped with rules for the initialization of and numerical operations on multi-complex numbers,  is introduced in this paper. In particular,  an $n$-complex library has been developed for C++, adhering to the  C++11 standards~\cite{Cppeleven}. In addition, a MATLAB version of this library has been constructed and will be made available along with this paper.
\subsubsection{The $n$-complex number class}
The  class of $n$-complex numbers  has been  developed under the core concept of class inheritance. In an arbitrary class defining an $n$-complex number, two private members \texttt{z1} and \texttt{z2} of the $(n-1)$-complex number class are included, per the recursive definition 
of (\ref{ncomplexdef}).  For example, a variable of the $bi$-complex data type ($\mathbb{C}^n$, $n=2$) consists of two private members \texttt{z1} and \texttt{z2} of the $single$-complex number data type (i.e. $z_1, z_2 \in \mathbb{C}$). 

For computations involving $n$-complex numbers (as those in the FDTD update equations), basic mathematical and linear algebra operators are overloaded, to accept $n$-complex numbers as input arguments. Since the $n$-complex class is defined recursively, these functions can be used repeatedly without duplicate declarations. The above characteristics of the  $n$-complex library allow users to define an arbitrary number of imaginary dimensions. In addition to basic operators, it is necessary to have new functions for the extraction of the real and multiple imaginary  parts of an $n$-complex variable. For example, functions \texttt{Im1(z), Im2(z), Im12(z)} return  $\text{Im}_1(z) \equiv y$, $\text{Im}_2(z) \equiv z$, 
$\text{Im}_{12}(z) \equiv \tau$ of a $bi-$complex variable $z=x + j_1 y + j_2 z + j_{12}\tau$, respectively.
\subsubsection{Using the $n$-complex number class in FDTD}
The $n$-complex number class is employed to implement MCSD in standard FDTD codes as follows. First, all design parameters
of interest are declared  as  $n$-complex, consisting of two  $(n-1)$-complex number arrays $z_1, z_2$,  with the same number of 
cells as that of the computational domain. Let us consider an example where second order field derivatives with respect to the 
dielectric permittivity of the substrate of a microstrip structure are computed. Then, the array of the relative dielectric permittivities of 
the ($i,m,k$) Yee cell is defined as: 
\begin{equation}
\epsilon_r^{i,m,k} = z_1^{i,m,k} + j_2 z_2^{i,m,k}, \: \: z_1,\: z_2 \in \mathbb{C}  
\end{equation}
where the two arrays $z_1$, $z_2$ are: 
\begin{equation}
 z_1^{i,m,k}  = \left\{ 
\begin{array}{ll}
\epsilon_r^{\mathrm{substrate}} + j_1 h, & (i,m,k) \: \text{in substrate} \\[.3cm]
1, & (i,m,k) \: \text{in air} 
\end{array} 
 \right.
\end{equation}
and 
\begin{equation}
 z_2^{i,m,k}  = \left\{ 
\begin{array}{ll}
h & (i,m,k) \: \text{in substrate} \\[.3cm]
0, & (i,m,k) \: \text{in air} 
\end{array} 
 \right.
\end{equation}
Since the $n$-complex number class is included in the programming environment, the electric and magnetic field components
are transformed to the $n$-complex number data type if they are associated with perturbed cells. This
self-defined $n$-complex class is compatible with any existing FDTD codes in both MATLAB and C++ environments. In the following, FDTD codes loaded with the $n$-complex number class will be referred to as MCSD-FDTD. In the next sections, MCSD-FDTD is evaluated in 
terms of its accuracy and applied to a 3-D simulation scenario.
\subsubsection{Defining the $n$-complex number class in matrix form}
The alternative implementation 
\textcolor{black}{
\subsection{Jacobian and Hessian Matrix Computation}
This section further studies the cost of computing Jacobian and Hessian matrices via CSD- and MCSD-FDTD.
\subsubsection{Jacobian matrix}Consider a 3-D FDTD simulation discretized by $I\times M \times K$ cells, and $I+M+K=P$. Then, at each time step, a total of $6P$ field components are computed, where $p=1, 2, \cdots P$ is a Yee cell index corresponding to the triad $(i,m,k)$. The first-order derivatives of these field components with respect to $N$ design variables $\bs{\xi} = \left[ \xi_1, \xi_2, \cdots \xi_N  \right]^T$ are written in the form of a $P \times N$ Jacobian matrix:
 \begin{equation}
\nabla_{\bs{\xi}}\tilde{\textbf{E}}^{x,y,z}, \nabla_{\bs{\xi}}\tilde{\textbf{H}}^{x,y,z}
 \end{equation}
where
 \begin{equation}
\nabla_{\bs{\xi}}\tilde{\textbf{E}}^x \equiv
 \begin{bmatrix}
 \dfrac{\partial {E^x_1}}{\partial \xi_1} & \dfrac{\partial {E^x_1}}{\partial \xi_1} & \dots  &  \dfrac{\partial {E^x_1}}{\partial \xi_N} 
\\[.4cm]
    \dfrac{\partial {E^x_2}}{\partial \xi_1} & \dfrac{\partial {E^x_2}}{\partial \xi_2} & \dots  &  \dfrac{\partial {E^x_2}}{\partial \xi_N} \\
\vdots & \vdots & \ddots & \vdots \\
\dfrac{\partial {E^x_P}}{\partial \xi_1} & \dfrac{\partial {E^x_P}}{\partial \xi_1} & \dots  &  \dfrac{\partial {E^x_P}}{\partial \xi_N} \\
\end{bmatrix}
 \end{equation}
To derive these matrices, a total of $N$ CSD-FDTD simulations are needed. In each independent CSD-FDTD simulation, one imaginary perturbation is assigned to a specific parameter $\xi_i$. In contrast, by using the CFD-FDTD method, at least $2N$ CFD-FDTD simulations are needed.
\subsubsection{Hessian matrix}
In addition to first-order derivatives, the Hessian matrix of a field component can be further derived by using MCSD-FDTD with bi-complex numbers. To obtain the Hessian matrix of $E^z_p$ with respect to $N$ design parameters $\bs{\xi} = \left[ \xi_1, \xi_2, \cdots \xi_N  \right]^T$, written as,
\begin{equation}
  \textbf{H} \left( E^z_p \right) = 
  \begin{bmatrix}
    \dfrac{{\partial}^2 E^z_p }{\partial {\xi_1}^2} & \dfrac{{\partial}^2 E^z_p }{\partial \xi_1 \partial \xi_2} & \dots & \dfrac{{\partial}^2 E^z_p }{\partial \xi_1 \partial \xi_N}   
 \\[.4cm]
    \dfrac{{\partial}^2 E^z_p }{\partial \xi_2 \partial \xi_1} & \dfrac{{\partial}^2 E^z_p }{\partial {\xi_2}^2 } & \dots & \dfrac{{\partial}^2 E^z_p }{\partial \xi_2 \partial \xi_N} 
\\
\vdots & \vdots & \ddots & \vdots
\\
\dfrac{{\partial}^2 E^z_p }{\partial \xi_N \partial \xi_1} & \dfrac{{\partial}^2 E^z_p }{\partial \xi_N \partial \xi_2 } & \dots & \dfrac{{\partial}^2 E^z_p }{\partial {\xi_N}^2} 
  \end{bmatrix},
\end{equation}
at least $4N^2$ CFD-FDTD simulations are needed. Alternatively, a total of $N^2$ MCSD-FDTD simulations can be performed, with two imaginary perturbations $j_1h_1$, $j_2h_1$ assigned to $\xi_i$, $\xi_j$ ($i,j \in [1,N]$) respectively in each independent run. 
}
\begin{table*}[]
\centering
\caption{Operation Count for The Computation of High-order Derivatives 
\protect\\(Add./Sub. : Additions/Subtractions; Mult. : Multiplications)}
\label{time_table}
\begin{threeparttable}
\begin{tabular}{|c|c|c|c|c|c|c|}
\hline
\multirow{2}{*}{} & \multicolumn{2}{c|}{CFD-FDTD\tnote{$\star$}} & \multicolumn{2}{c|}{MCSD-FDTD} & \multicolumn{2}{c|}{Iterative CSD-FDTD} \\ \cline{2-7} 
                  & Add./Sub.       & Mult.       & Add./Sub.        & Mult.       & Add./Sub.            & Mult.            \\ \hline
\begin{tabular}[c]{@{}l@{}}\\$\dfrac{{\partial}^n \tilde{\textbf{E}}}{\partial \xi^n}$, $1\leqslant n\leqslant N$\\ $ $\end{tabular}      & \multicolumn{2}{c|}{$2N$ }    & $2^N$    & 3$^N$ & \multicolumn{2}{c|}{$2N$}         \\ \hline
\begin{tabular}[c]{@{}l@{}}\\$\dfrac{{\partial}^{m_1+m_2+\dots+m_N} \tilde{\textbf{E}}}{\partial^{m_1} {{\xi}_1}\partial^{m_2} {{\xi}_2} \dots \partial^{m_N} {{\xi}_N}}$, $m_i=0,1$\\$ $\end{tabular}    & \multicolumn{2}{c|}{$3^N$}    & $2^N$  & $3^N$  & \multicolumn{2}{c|}{$2N$}  \\ \hline
\end{tabular}
\begin{tablenotes}
     \item[$\star$] Excluding FDTD runs to ensure the convergence of the result.
\end{tablenotes}  
\end{threeparttable}
\end{table*}
\textcolor{black}{
\subsection{High-order Partial Mixed Derivatives}
In addition to computing the second-order partial mixed field derivatives in the Hessian matrix, MCSD-FDTD is capable of producing $N$-order mixed partial derivatives of field components with respect to $N$ design variables $\bs{\xi} = \left[ \xi_1, \xi_2, \cdots \xi_N  \right]^T$:
\begin{equation}
\dfrac{{\partial}^N E^z_p}{\partial \xi_1 \partial \xi_2 \dots \partial \xi_N },
\label{high-order mixed1}
\end{equation}
$N$ imaginary perturbations $j_1h_1, j_2h_2, \dots, j_Nh_N$ are assigned to each variable respectively. In addition to $N$-th order mixed partial derivative, $N'$-th order partial mixed derivatives with respect to $\bs{\xi}' = \left[ \xi_1, \xi_2, \cdots \xi_{N'}  \right]^T$, where $N' = 1, \cdots N-1$, are computed at the same time in this MCSD-FDTD simulation. If CFD-FDTD is used to find all these derivatives, at least $3^N$ real-valued FDTD simulations are needed.
\subsection{High-order Derivatives with an Iterative CSD-FDTD Method}
The computational redundancy in finding high-order field derivatives with respect to one variable is found by examining (\ref{mcsdoutput}), in which both $\text{Im}_{1} \left\{ \tilde{\textbf{E}}_{i,m,k}^n \: \right\}/h_1$ and  $\text{Im}_{2} \left\{ \tilde{\textbf{E}}_{i,m,k}^n \: \right\}/h_2$ approximate $\partial \textbf{E}_{i,m,k}^n / \partial \epsilon_r$. The additional computation overhead becomes heavier if more imaginary perturbations are assigned to one variable. For example, to find $\partial^3 \textbf{E}_{i,m,k}^n / \partial {\epsilon_r}^3$, three imaginary perturbations $j_1h_1, j_2h_1, j_3h_1$are assigned to $\epsilon_r$. The MCSD approximation to the third-order derivative is simply:
\begin{equation}
\dfrac{{\partial}^3 \left\{ \tilde{\textbf{E}}_{i,m,k}^n \: \right\}}{\partial {\epsilon_r}^3} \approx \dfrac{\text{Im}_{123} \left\{ \tilde{\textbf{E}}_{i,m,k}^n \: \right\}}{{h_1}^3}
\end{equation}
In addition to the third-order derivative, multiple approximations to first and second-order derivatives are redundantly computed:
\begin{equation}
\begin{array}{l}
\text{Im}_{1,2,3} \left\{ \tilde{\textbf{E}}_{i,m,k}^n \: \right\}/h_1 
\approx \partial \textbf{E}_{i,m,k}^n / \partial \epsilon_r \\[.2cm]
\text{Im}_{12,23,13} \left\{ \tilde{\textbf{E}}_{i,m,k}^n \: \right\}/{h_1}^2 
\approx {\partial}^2 \textbf{E}_{i,m,k}^n / \partial {\epsilon_r}^2 \\[.2cm]
\end{array}
\label{high-order redundant}
\end{equation}
}
\textcolor{black}{
To alleviate the computational overhead of MCSD-FDTD, an alternative approach using only one imaginary perturbation per design parameter is proposed. Based on that, derivatives up to any order can be computed iteratively, with a ``marching-in-order" approach. This iterative CSD-FDTD method is elucidated through the update equation of the $x$-component of the electric field:
\begin{equation}
E_{i',m,k}^{x,n+1} =  E_{i',m,k}^{x,n}  + \mathcal{J},
\label{iterativeCSD-1}
\end{equation}
where
\begin{multline}
\mathcal{J} = \dfrac{\Delta t}{\epsilon_0 \epsilon_r } \times
 \Bigg( \dfrac{H_{i',m',k}^{z,n'}- H_{i',m'-1,k}^{z,n'}}{\Delta y} 
 -  \dfrac{H_{i',m,k'}^{y,n'}-
H_{i',m,k'-1}^{y,n'}}{\Delta z} \Bigg)
\end{multline}
Direct differentiation of (\ref{iterativeCSD-1}) with respect to $\epsilon_r$ yields:
\begin{equation}
\dfrac{ {\partial} E_{i',m,k}^{x,n+1} }{ \partial {\epsilon_r} } =  \dfrac{ {\partial} E_{i',m,k}^{x,n} }{\partial {\epsilon_r}} + \mathcal{K},
\label{iterativeCSD-2}
\end{equation}
where 
\begin{multline} 
\mathcal{K}=
\\[.4cm]
\dfrac{ \Delta t}{\epsilon_0 \epsilon_r}
\Bigg( \dfrac { \dfrac{ {\partial} H_{i',m',k}^{z,n'} }  {\partial  {\epsilon_r} } - 
\dfrac{{\partial} H_{i',m'-1,k}^{z,n'} } { \partial  {\epsilon_r} } } {\Delta y}
 -  \dfrac{ \dfrac{ {\partial} H_{i',m,k'}^{y,n'} } {\partial  {\epsilon_r}} - 
 \dfrac{{\partial} H_{i',m,k'-1}^{y,n'} } { \partial  {\epsilon_r} } } {\Delta z} \Bigg)
\\[.4cm]
 - \dfrac{\Delta t}{ \epsilon_0 {\epsilon_r}^2}
 \Bigg(\dfrac{ H_{i',m',k}^{z,n'} - H_{i',m'-1,k}^{z,n'}}{\Delta y} 
- \dfrac{ H_{i',m,k'}^{y,n'} - H_{i',m,k'-1}^{y,n'}}{\Delta z} \Bigg)
\\[.4cm]
\end{multline}
Introducing an  imaginary perturbation $j_1h_1$ to $\epsilon_r$ in (\ref{iterativeCSD-1}) and (\ref{iterativeCSD-2}) 
turns them into complex-valued equations. Then,  (\ref{iterativeCSD-2}) is updated iteratively by the real and imaginary parts of 
$\mathcal{J}$ from (\ref{iterativeCSD-1}), as: 
\begin{multline}
\dfrac{ {\partial} E_{i',m,k}^{x,n+1} }{ \partial {\epsilon_r} } =  \dfrac{ {\partial} E_{i',m,k}^{x,n} }{\partial {\epsilon_r}} + 
\\[.4cm]
\dfrac{ \Delta t}{\epsilon_0 \left( \epsilon_r + j_1h_1 \right)} \times \text{Re} \left\{ 
\mathcal{J} \right\} - \dfrac{ \Delta t}{\epsilon_0 \left( \epsilon_r + j_1h_1 \right)^2} \times \text{Im}_1\left\{ \mathcal{J} \right\}
\end{multline}
Higher order field derivatives can be further derived following the same scheme. In general, the added term ${\mathcal{K}}$ that appears in the $p$-th order field derivative is:
\begin{multline}
\dfrac{\Delta t}{\epsilon_0}\sum_{q=0}^{p-1} \frac{p!}{(p-q)!q!}
\\
\times \Bigg[
\dfrac{{\partial}^{p-q}}{\partial {\epsilon_r}^{p-q} \Delta y}
\Bigg( \dfrac{ {\partial}^q H_{i',m',k}^{z,n'} }  {\partial  {\epsilon_r}^q } - 
\dfrac{{\partial}^q H_{i',m'-1,k}^{z,n'} } { \partial  {\epsilon_r}^q }
\Bigg)-
\\
\dfrac{{\partial}^{p-q}}{\partial {\epsilon_r}^{p-q} \Delta z}
\Bigg( \dfrac{ {\partial}^q H_{i',m,k'}^{y,n'} } {\partial  {\epsilon_r}^q} - 
 \dfrac{{\partial}^q H_{i',m,k'-1}^{y,n'} } { \partial  {\epsilon_r}^q  } \Bigg) \Bigg]
\end{multline}
Therefore, to find the $N-$th order field derivative with respect to a particular variable, $N$ additional terms are needed in ${\mathcal{K}}$, and these terms are available from the update equations of the $(N-1)$-th order field derivatives. 
}

\textcolor{black}{
Finally, 
the operation count for MCSD-FDTD is summarized and compared to that of CFD-FDTD, relative to a single real-valued FDTD simulation, in Table \ref{time_table}. Parallelization or multi-threading, both viable options to accelerate these computations,  are not considered in this analysis. For one MCSD-FDTD simulation with field arrays in $\mathbb{C}^n$, 
there are  $2^N$ times the additions/subtractions and  $3^N$ times the multiplications of  a real-valued FDTD simulation. 
}
\section{MCSD-FDTD: Validation}
\subsection{The relation between complex step and accuracy}
\begin{figure}[]
	\centering
 		\subfigure[]{\includegraphics[width=8cm]{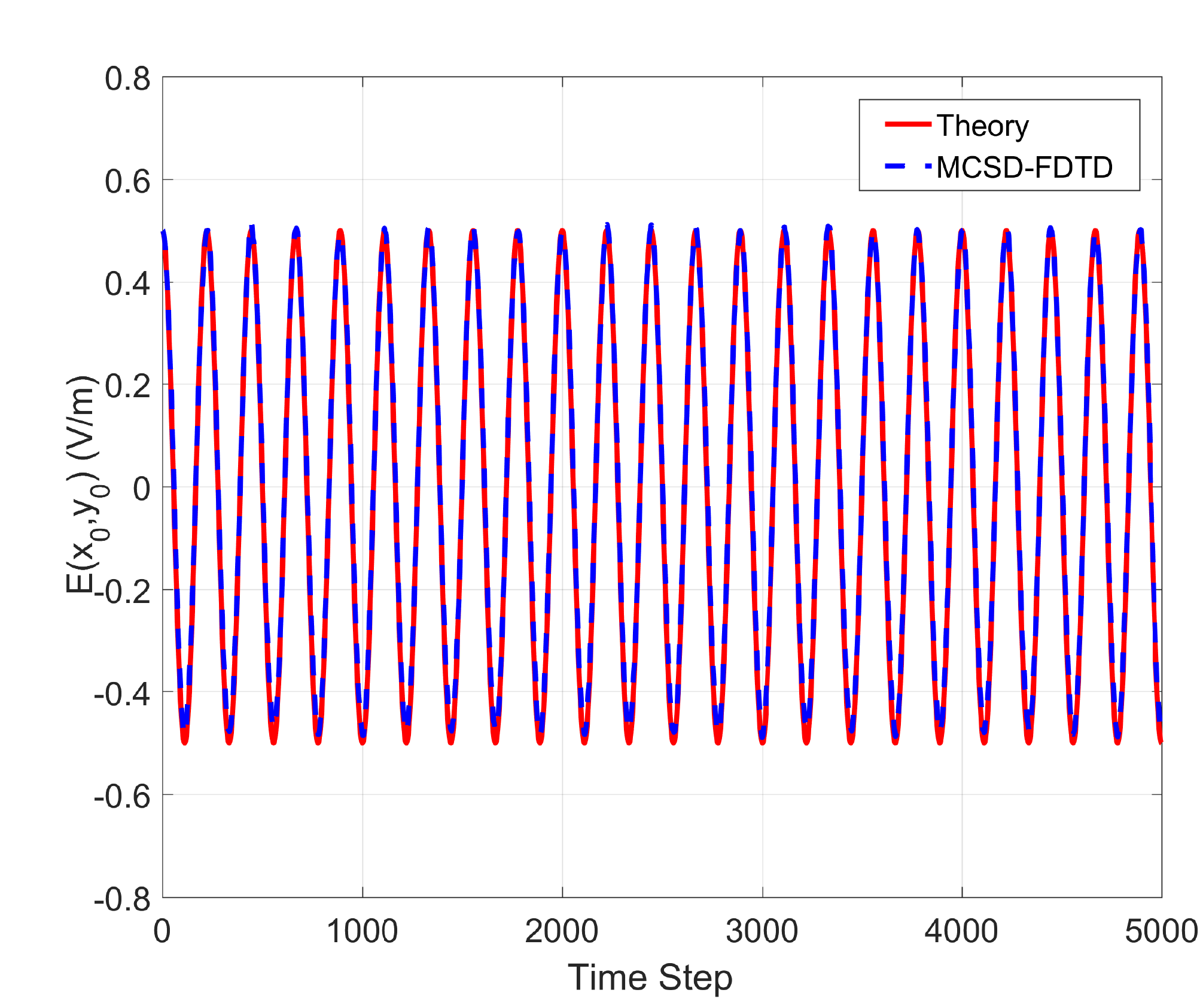}}
 		\subfigure[]{\includegraphics[width=8cm]{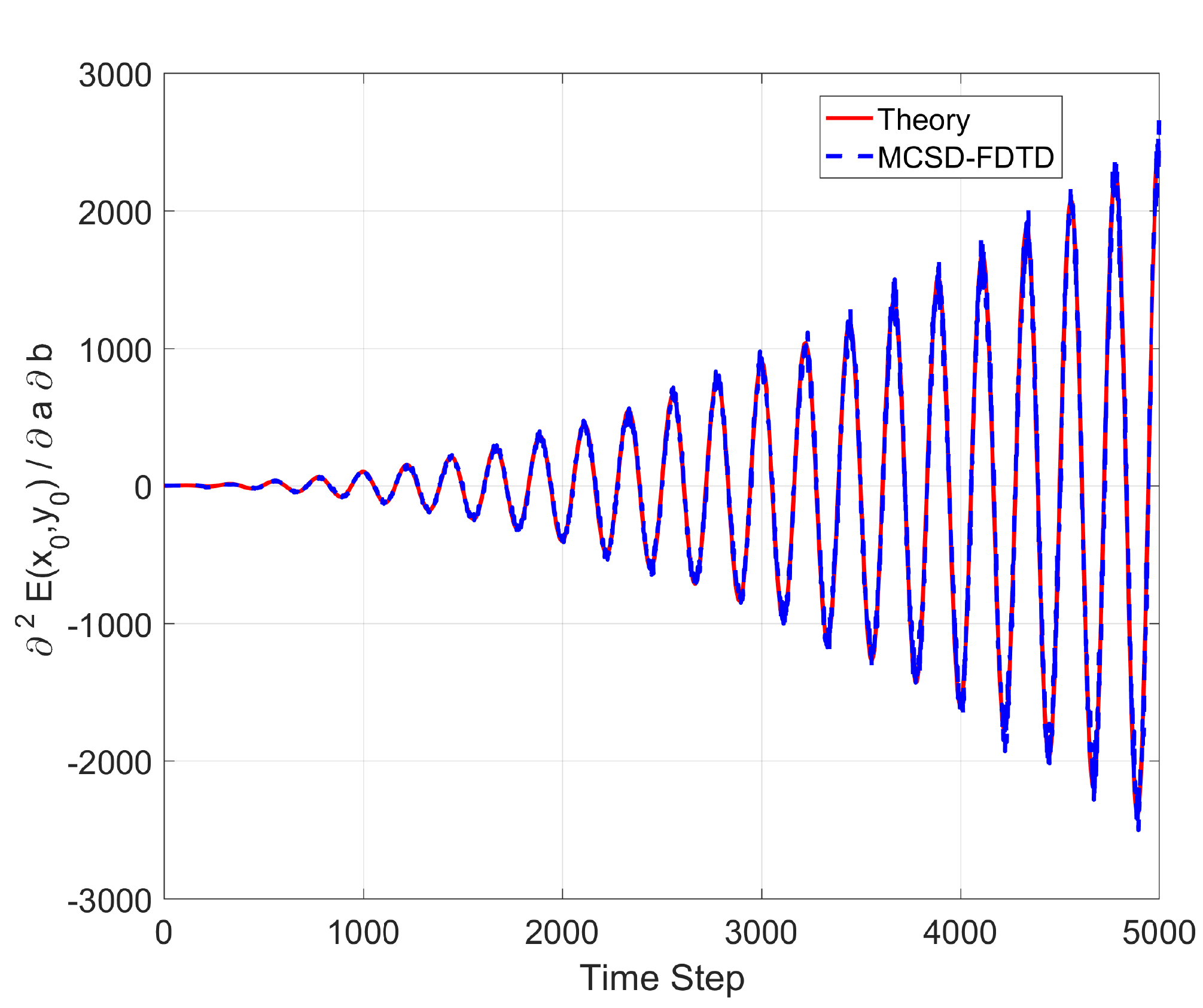}}
		\caption{Comparison of analytical solution and MCSD-FDTD for (a) the electric field and (b) its second-order derivative 
		with respect to the dimensions of a rectangular metallic cavity,  at a sampling point within the cavity, 
		in the time-domain. }
		\label{cavity_timedomain}
\end{figure}
The accuracy of MCSD-FDTD is first tested in a two-dimensional case study, where field derivatives are analytically available: the  derivatives of electromagnetic field components within a rectangular, air-filled, metallic cavity, with respect to the dimensions of the  cavity. The cavity is discretized by a uniform mesh of $150\times 100$  Yee cells with $\Delta x = \Delta y=\Delta = 1$ mm. Transverse  electric (TE) modes with ($H_x$, $H_y$, $E_z$) are considered. The electric field of the ($m$, $n$) mode is: 
\begin{equation}
E^{(m,n)}_z (x, y ,t ) = E_0  \sin{\frac{m \pi x}{a}} \sin{\frac{n \pi  y}{b}}  \cos{\left( 2 \pi f_{m,n} t \right)}
\end{equation}
with 
\begin{equation}
f_{m,n} = \dfrac{1}{2 \sqrt{\epsilon_0 \mu_0}} \sqrt{ \left( \dfrac{m}{\alpha}\right)^2 + \left( \dfrac{n}{b}\right)^2 }
\end{equation}
In the expressions above, $\alpha=15$ cm is the width and $b=10$\,cm is the height of the cavity. Based on these, 
the derivatives of the ($m,n$) modal fields with respect to $\alpha$ and $b$ are found and compared to their 
numerically computed values, via MCSD-FDTD, in the time-domain. To that end, the electric field of the
($m,n$) mode is injected as an initial condition (at $t=0$) and its time evolution is simulated.


In this example, first and second order derivatives of the field components with respect to $\alpha$ and $b$ are 
computed. Therefore, the cells at the edges of the cavity (in both the $x$- and $y$-direction) are perturbed by 
complex steps $j_1 h_1$ and $j_2 h_2$, and set to $\Delta x' \equiv  \Delta x \left( 1 + j_1 h_1 \right)$, 
$\Delta y' \equiv  \Delta y \left(1  + j_2 h_2 \right)$, with $h_1=h_2=10^{-5}$.  
%
The corresponding electric field update equation is: 
\begin{multline}
E_{i',m,k}^{z,n+1} =  E_{i',m,k}^{z,n}  +  \dfrac{\Delta t}{\epsilon_r  \Delta y'}  \left(
H_{i',m',k}^{x,n'}- H_{i',m'-1,k}^{x,n'} \right) \\[.3cm]
  -   \dfrac{\Delta t}{\epsilon_r  \Delta x'} \left(  H_{i',m,k'}^{z,n'} - 
H_{i',m,k'-1}^{z,n'} \right)
\end{multline} 
The results of MCSD-FDTD are compared to the analytical ones in Fig. \ref{cavity_timedomain}, for the electric 
field of the ($1,1$) cavity mode. The real part of the electric field in MCSD-FDTD accurately reproduces the electric field of the unperturbed problem (sampled at the center of the cavity), as shown in Fig. \ref{cavity_timedomain}(a). Moreover,  Fig. \ref{cavity_timedomain}(b) shows the second order derivative
of the electric field at the same sampling point, with respect to $\alpha$ and $b$, found by MCSD-FDTD, along with its 
analytical expression. Note that this derivative  grows with time, as: 
\begin{subequations}
\begin{align}
\dfrac{\partial \cos{(2 \pi f_{mn} t )}}{\partial \alpha} = -\sin{(2 \pi f_{mn} t )} 2 \pi t \dfrac{\partial f_{m,n}}{
\partial \alpha} \\
\dfrac{\partial \cos{(2 \pi f_{mn} t )}}{\partial b} = -\sin{(2 \pi f_{mn} t )} 2 \pi t \dfrac{\partial f_{m,n}}{
\partial b} 
\end{align}
\end{subequations}

The temporal evolution of $\partial^2 E_z / \partial \alpha \partial b$ is accurately represented by the MCSD-FDTD 
solution, in agreement with its analytical counterpart. To further elucidate the accuracy of the MCSD-based derivative 
calculation and to compare it to the conventional alternative of centered finite-differences (CFD), the analytical expression of  $\partial^2 E_z / \partial \alpha \partial b$ is employed to characterize the accuracy of 
MCSD and CFD with respect to the step size $h$. The following $l_{\infty}$ error norm is used: 
\begin{equation}
\max_{n} \max_{i,j} \dfrac{| E_{i,j}^{z,n}(\text{MCSD/CFD}) - E_{i,j}^{z,n}(\text{theory})|}{|E_{i,j}^{z,n}(\text{theory})|} 
\label{errornorm}
\end{equation}
This registers the maximum error over space and time within 5,000 time steps of the simulation. 
\begin{figure*}[]
	\centering
		\subfigure[]{\includegraphics[width=8.2cm]{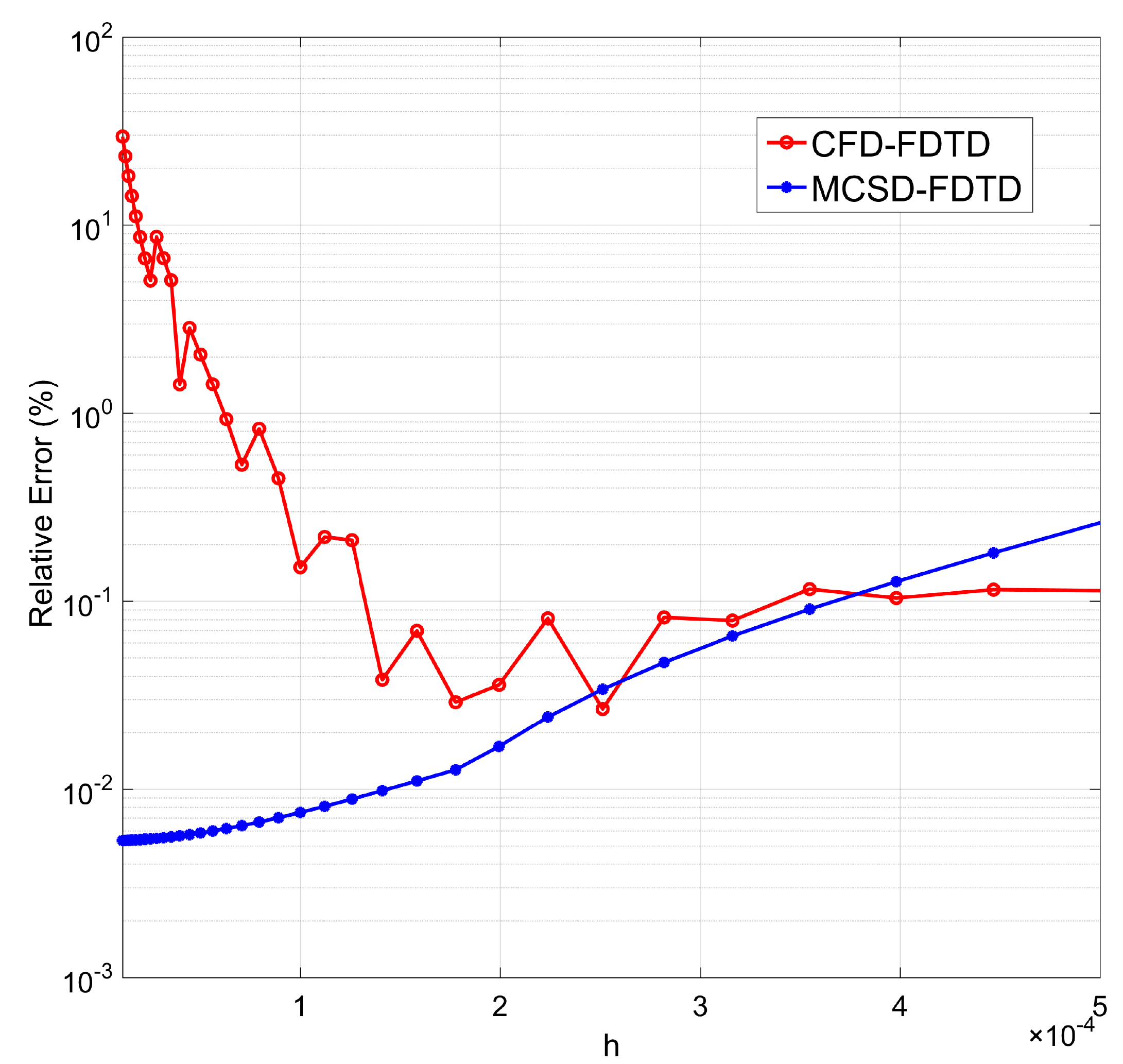}
		\label{cavity_error_compare}}
    	\subfigure[]{
    	\includegraphics[width=8cm]{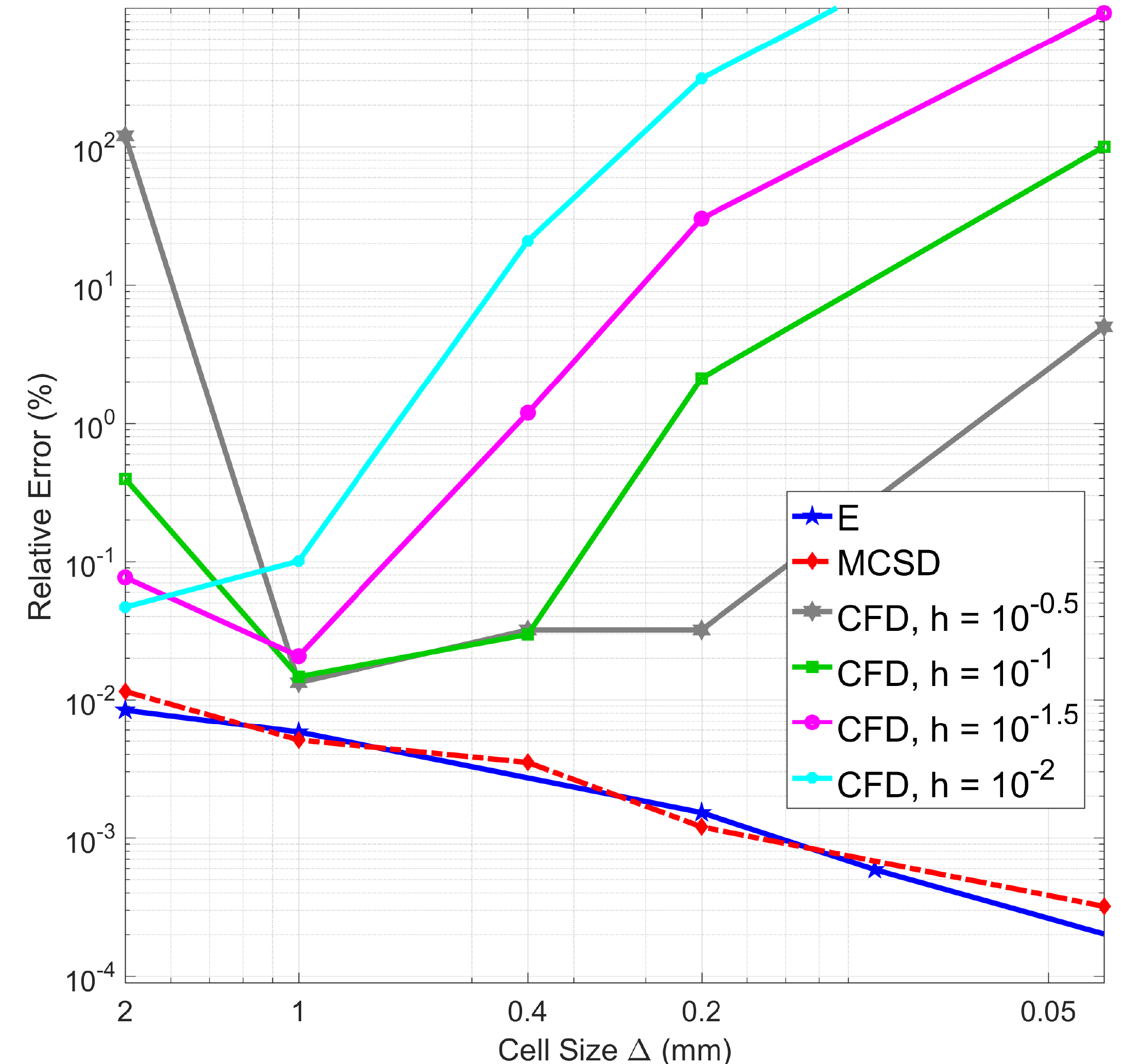}
    	\label{cavity_error_compare_different_mesh}}
\caption{
\textcolor{black}{
Relative error, according to (\ref{errornorm}), of the numerical sensitivity $\partial^2 E / \partial a \partial b $ of the electric field inside a rectangular cavity, computed via FDTD with central finite differences (CFD), as well as the multi-complex step derivative (MCSD) approximation. In (a), the Yee cell size $\Delta$ is set at $1$ mm and step-size $h$ varies from $5\times10^{-4}$ to $10^{-5}$. In (b), the cell size varies from $2$ to $0.4$ mm, and the relative error of the electric field $E$ computed by real-valued FDTD is included.}}
\end{figure*}

\begin{figure}[]
	 \centering
		\subfigure[]{\includegraphics[width=8cm]{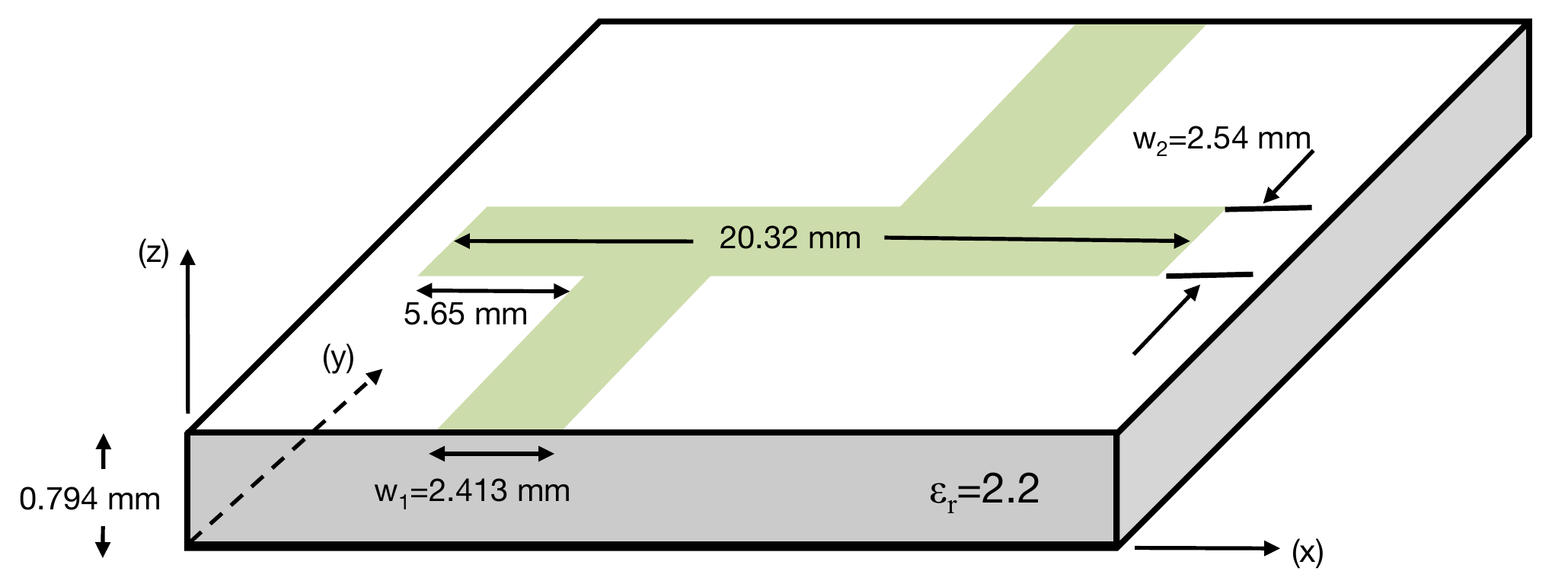} \label{filtergeom}}
		 \subfigure[]{\includegraphics[width=8cm]{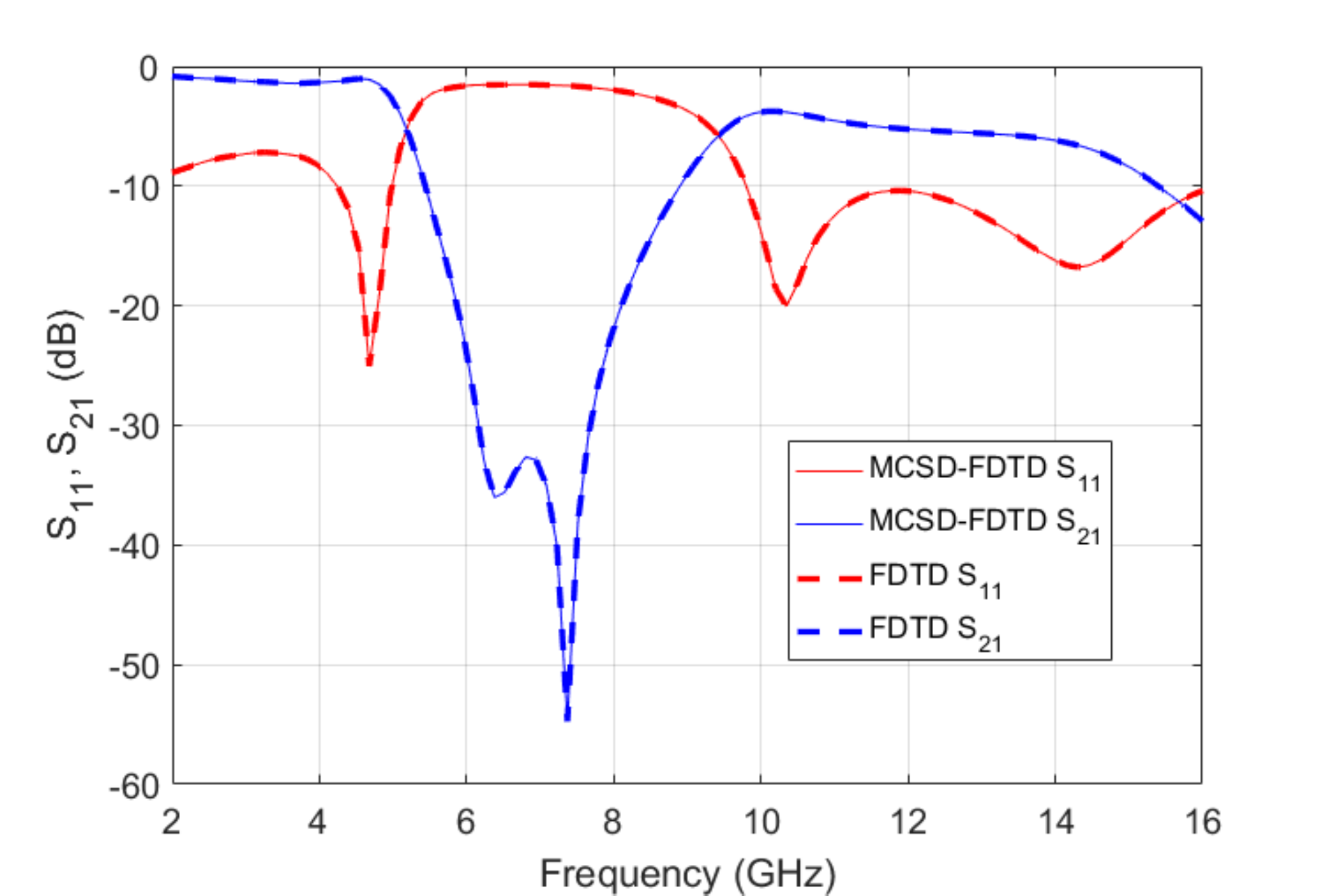} \label{filterS11S21}}
		\caption{In (a), the geometry of a microstrip filter, from \cite{sheen3DFDTD}, is shown. The $S$-parameters of the filter, computed via 
		standard (unperturbed) FDTD and MCSD-FDTD are shown in (b).}
		\label{filter3rd}
\end{figure}
\begin{figure*}
	\centering
		\subfigure[]{\includegraphics[width=7.8cm]{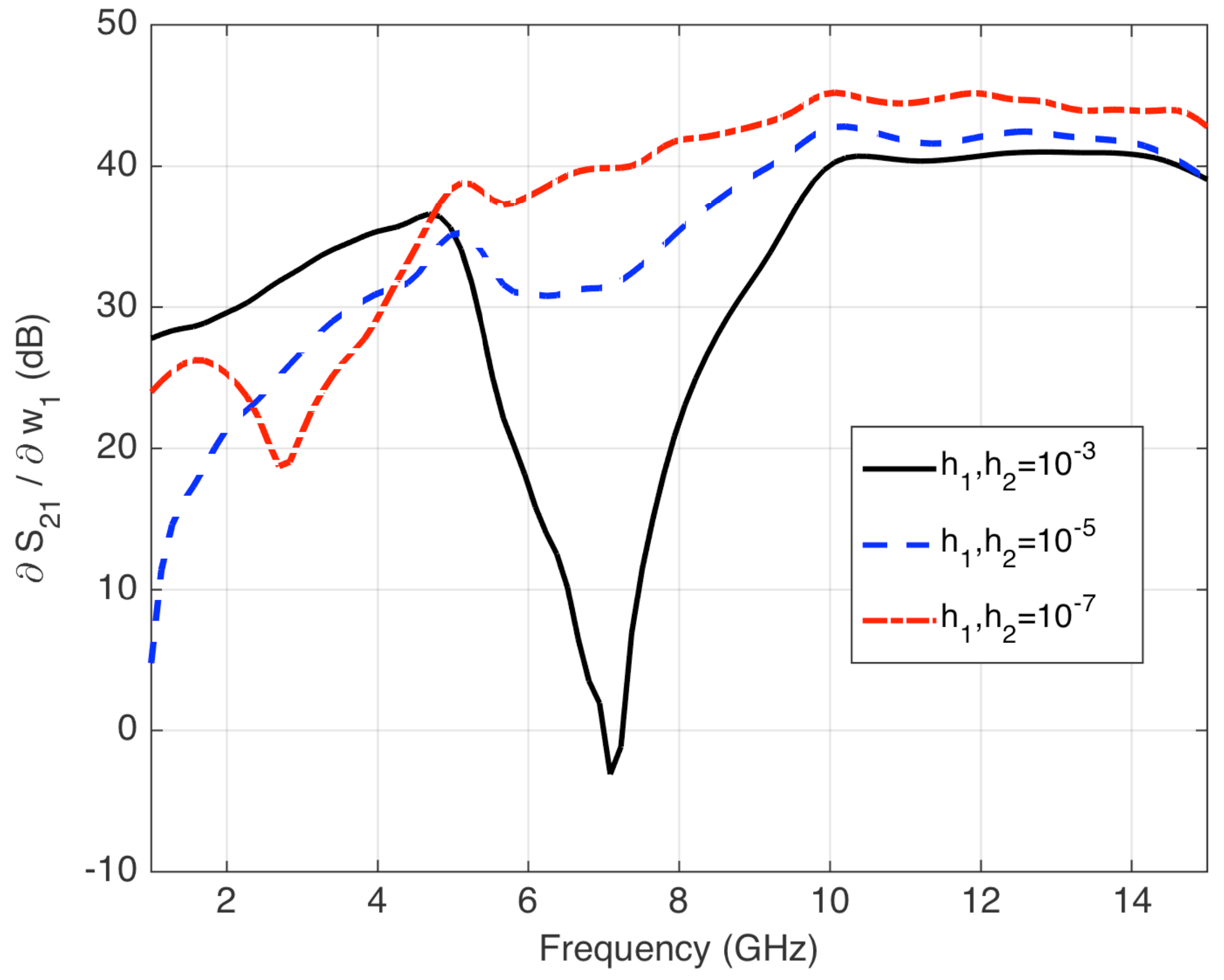}}
		 \subfigure[]{\includegraphics[width=8cm]{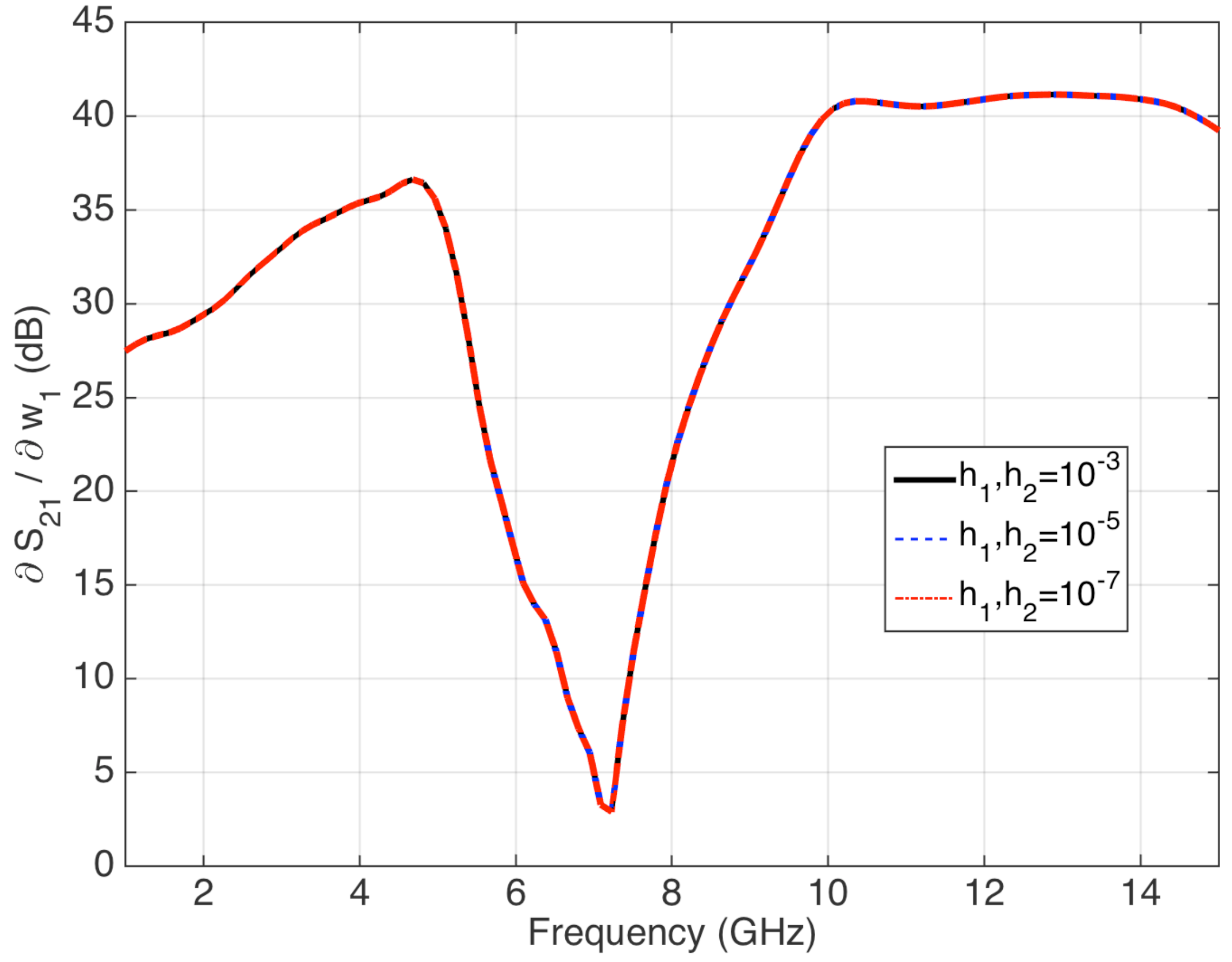}}
		 \subfigure[]{\includegraphics[width=8cm]{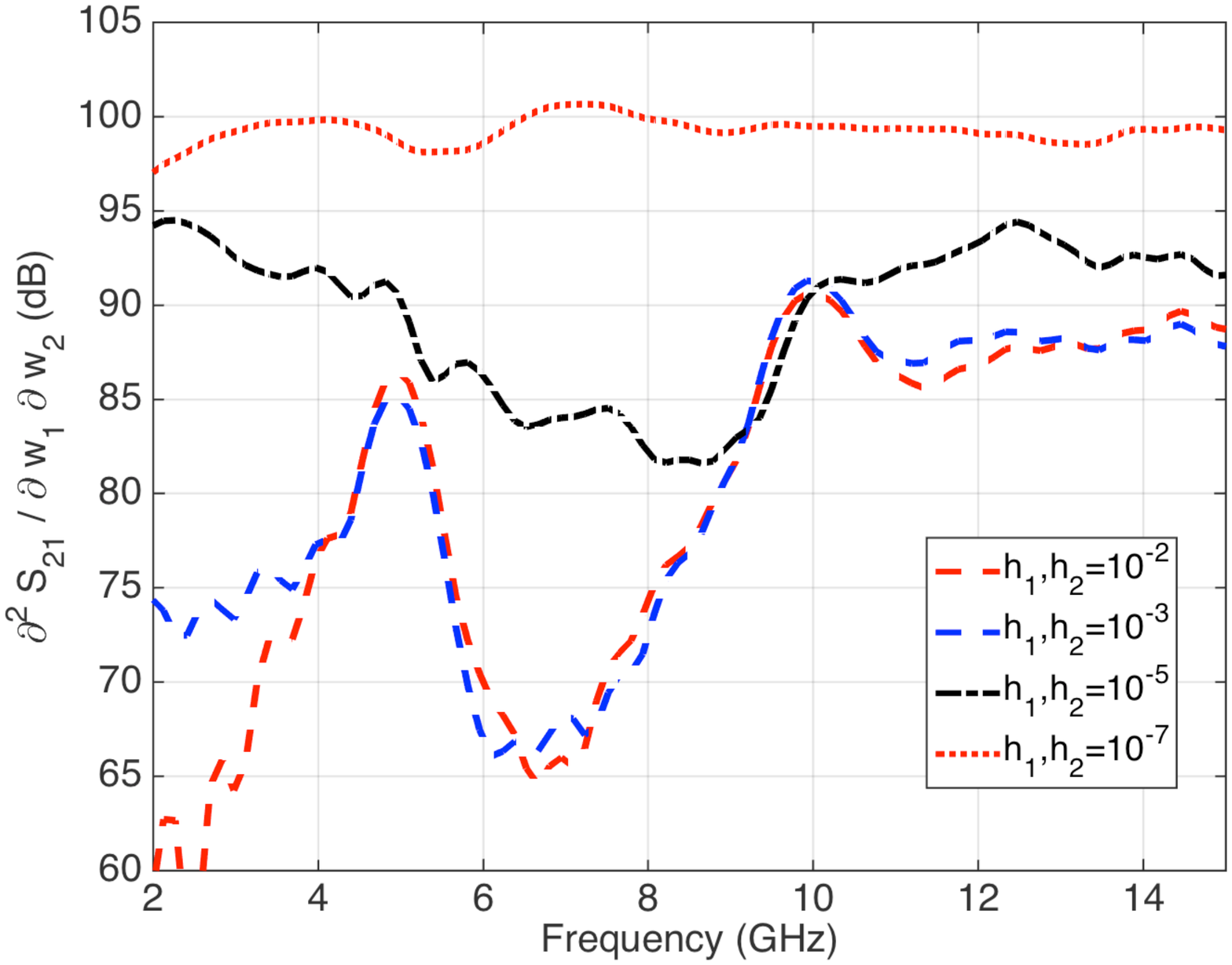}}
		 \subfigure[]{\includegraphics[width=8cm]{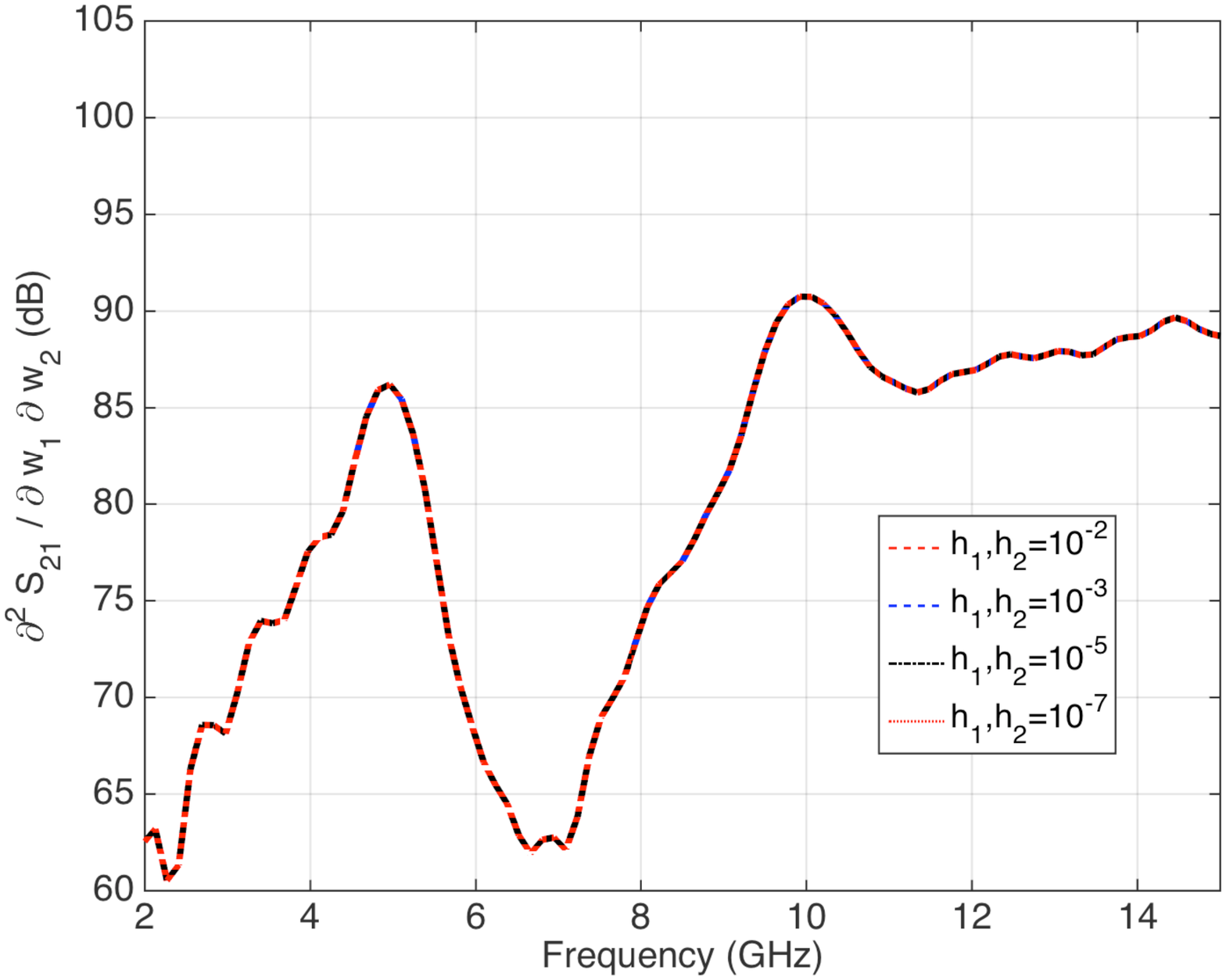}}
		\caption{Numerical sensitivities of the microstrip low-pass filter: $\partial S_{21} / \partial w_1$ (computed via FDTD combined with CFD, in (a),  and MCSD, in (b)), and $\partial^2 S_{21} / \partial w_1 \partial w_2$ (computed via FDTD combined with CFD, in (c), and MCSD, in (d)).}
		\label{filtersens}
\end{figure*}
In Fig.~\ref{cavity_error_compare}, this relative error norm (\ref{errornorm}) is plotted with respect to the step size $h$. 
Notably, the general error trends that were found in the theoretical analysis of CFD and MCSD in section II
appear here as well. The important difference is that these errors are now super-imposed to the FDTD dispersion 
errors. Hence, MCSD-FDTD cannot reach machine accuracy levels. Yet, it does outperform CFD by an order of magnitude, 
in terms of accuracy. It is also evident that CFD fails to converge, as $h$ is reduced; $h=10^{-4}$ appears to be a
breakpoint, where subtraction errors start to dominate. On the other hand, MCSD has a stable performance. As derived in (\ref{CSD2}) and (\ref{CFD2}), the leading error term of second-order derivative is proportional to $h^2/3$ and $h^2/12$ by MCSD and by CFD respectively. Hence, the MCSD-FDTD is prone to higher relative error than CFD-FDTD when $h>10^{-4}$. As MCSD-FDTD produces converged results as step-size is reduced, it is
allowed to use $h=10^{-5}$ in our computations in Fig. \ref{cavity_timedomain}, accurately recovering the unperturbed solution through the real part of the fields. 

In terms of execution time, one run of MCSD-FDTD took 16.4 secs. On the other hand, nine FDTD simulations, each taking 2.4 secs are needed 
to compute the same derivatives as MCSD-FDTD via CFD. Hence, the total execution time was  21.6 secs, {\em excluding} runs that would be 
needed to choose the step $h$ used in the centered finite differences. 
\textcolor{black}{
\subsection{FDTD numerical dispersion and MCSD-FDTD}
The relation of the numerical error of field derivatives and FDTD dispersion error is further studied in this section. The Yee cell size, $\Delta$, of the uniform mesh of this cavity is varied from $0.4$ to $2$ mm. In Fig. \ref{cavity_error_compare_different_mesh}, the relative error of the electric field calculated by FDTD compared to the analytical solution decreases with mesh refinement, along with the FDTD dispersion error. Notably, the relative error of the second-order mixed partial field derivative computed by MCSD-FDTD also decreases quadratically. In contrast, the subtractive cancellation error associated with the CFD method constrains its accuracy even as the FDTD dispersion error is 
minimized. 
}
\textcolor{black}{
Based on Fig. \ref{cavity_error_compare} and Fig. \ref{cavity_error_compare_different_mesh}, the advantage of MCSD-FDTD over CFD-FDTD for the computation of field derivatives is substantial. MCSD-FDTD provides guaranteed field derivatives with second-order accuracy, and it can be improved by reducing the step-size or the Yee cell size in FDTD. Hence, there is no need for additional simulations just to assess the accuracy of the derivative approximation as in CFD-FDTD. Note also that MCSD-FDTD with a Yee cell size $\Delta = 2$mm has
the same  accuracy as CFD-FDTD with cell size $\Delta = 1$ mm, in terms of second-order field derivatives. 
}
\section{The 3-D MCSD-FDTD: Microstrip filter}
A microstrip filter geometry, originally studied in \cite{sheen3DFDTD} and reproduced here in Fig. \ref{filtergeom}, is now used to demonstrate a three-dimensional application of  MCSD-FDTD and its distinct advantages over standard finite-difference methods. In particular, we focus on sensitivities 
of the $S$-parameters of the filter, with respect to the widths $w_{1,2}$, shown in Fig. \ref{filtergeom}. 
To that end, the derivatives of the $S$-parameters, are expressed in terms of the derivatives of the scattered and transmitted fields at the ports of the filter with respect to the design 
parameters. In turn, these derivatives are found applying MCSD-FDTD. The relevant computations are further explained in the following. 

The structure is discretized with a mesh of $80\times100\times16$ cells, with $ \Delta x=0.4064$ mm, $\Delta y=0.4233$ mm,
$\Delta z = 0.265$mm. $4000$ time steps ($\Delta t = 0.441 $ps) are used for the extraction of $S$- parameters in the frequency domain. A Gaussian pulse with half-width $T = 15$ ps is used as a source excitation.
For the MCSD-FDTD based computation of the sensitivities with respect to $w_{1,2}$, the method of section \ref{geom-MCSD} is applied. For 
$w_1$, the $\Delta x$ of the cells with the largest $x$-index modeling the segment of width $w_1$ are perturbed by $j_1 h_1\Delta x$. Likewise, 
the $\Delta y$ of the cells with the largest $y$-index modeling the segment of width $w_2$ are perturbed by $j_2 h_2 \Delta y$. 
%

In particular, the Yee cell sizes in the $x$- and $y$-directions,  $\Delta x $ and $\Delta y$ respectively,   are
initialized as three-dimensional $80\times 100\times 16$ real-valued arrays. The $\Delta x $ sub-array corresponding to the edge of
the first microstrip along the $y$- direction is perturbed by  $j_1h_1$, namely $\Delta x(36,1:46,4)\equiv \Delta x(1 +j_1h_1)$.
The $\Delta x$ of the neighbouring cells become $\Delta x(37,1:46,4)\equiv \Delta x(1 - j_1h_1)$.
Similarly, the sub-array $\Delta y$ corresponding to the edge of middle patch in $x$- direction is modified as:  
$\Delta y(16:66,50,4)\equiv \Delta y(1 +j_2h_2)$  and $\Delta y(16:66,51,4)\equiv \Delta y(1 -j_2h_2) $, 
to implement the perturbation in $w_2$. No perturbations are needed in the $z-$ direction. Hence, $\Delta z$ is real everywhere. With this choice, 
MCSD-FDTD can provide the  field solution to the unperturbed problem, through the real part of the fields. This is shown in Fig. 
\ref{filterS11S21}, which demonstrates the excellent agreement in the $S_{11}$ and $S_{21}$ found via standard FDTD applied to the unperturbed problem  and MCSD-FDTD. 

Furthermore, Fig. \ref{filtersens} shows first and second order derivatives of $S_{21}$, computed via MCSD-FDTD and CFD, for values 
of $h_1$, $h_2$ that vary between $10^{-7}$ and $10^{-3}$. These results indicate the significant convergence problem of CFD and the 
robustness of MCSD. The gap between the two becomes even more significant in the case of $\partial ^2 S_{21} / \partial w_1 \partial w_2$, where 
the CFD results diverge before reaching any satisfactory level of convergence over the simulated frequency bandwidth. On the other hand, 
the MCSD-FDTD results remain practically unchanged as $h_1$ and $h_2$ vary. 

\textcolor{black}
{With regards to execution time, one run of MCSD-FDTD took 142 seconds for this 3-D microstrip filter simulation. Alternatively, nine FDTD simulations, each taking 13.3 seconds are needed to compute the same derivatives as MCSD-FDTD via CFD. The total execution time is therefore $13.3 \times 9 = 119.7 $ seconds, excluding the number of additional runs needed to ensure the convergence of CFD.}

\textcolor{black}{
The capability of MCSD-FDTD to compute field  derivatives up to any order is further presented through the following example. We assign one more imaginary perturbation to the permittivity sub-array corresponding to the substrate of the filter, that is, $\epsilon_r (1:80,1:100,1:3)\equiv \epsilon_r(1 +j_3h_3)$. The third-order derivatives of $S$-parameters with respect to widths $w_{1,2}$ and substrate permittivity $\epsilon_r$ are found following the generalized approximation in (\ref{MCSD}) and shown in Fig. \ref{filter3rd} (a). On the other hand, these three imaginary perturbations can be assigned to a particular parameter for the computation of high-order derivatives with respect to one variable. Here,  $\Delta z(1:80,1:100,3)\equiv \Delta z(1 +j_1h_1+j_2h_2+j_3h_3)$ and the neighbouring cells $\Delta x(1:80,1:100,4)\equiv \Delta z(1 - j_1h_1-j_2h_2-j_3h_3)$. Third-order derivatives of $S$-parameters with respect to substrate thickness $d$ are calculated and shown in Fig. \ref{filter3rd} (b). Notably, the computation of this derivative through CFD-FDTD was not possible. The current state-of-art in AVM-FDTD has not gone past second-order derivatives.}

\begin{figure}[]
	 \centering
	 \vspace{-1cm}
		\subfigure[]{\includegraphics[width=8cm]{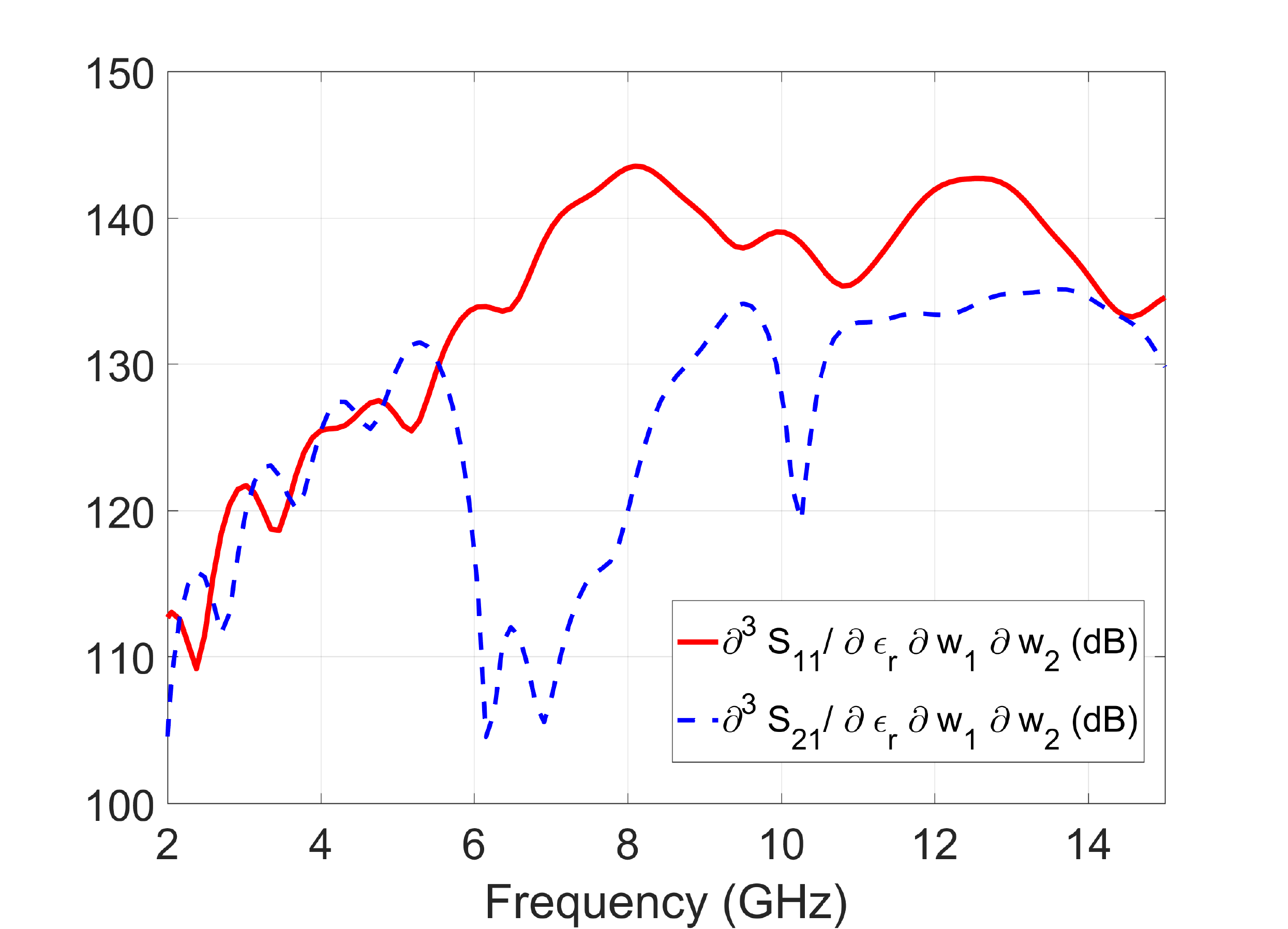}}
		 \subfigure[]{\includegraphics[width=8cm]{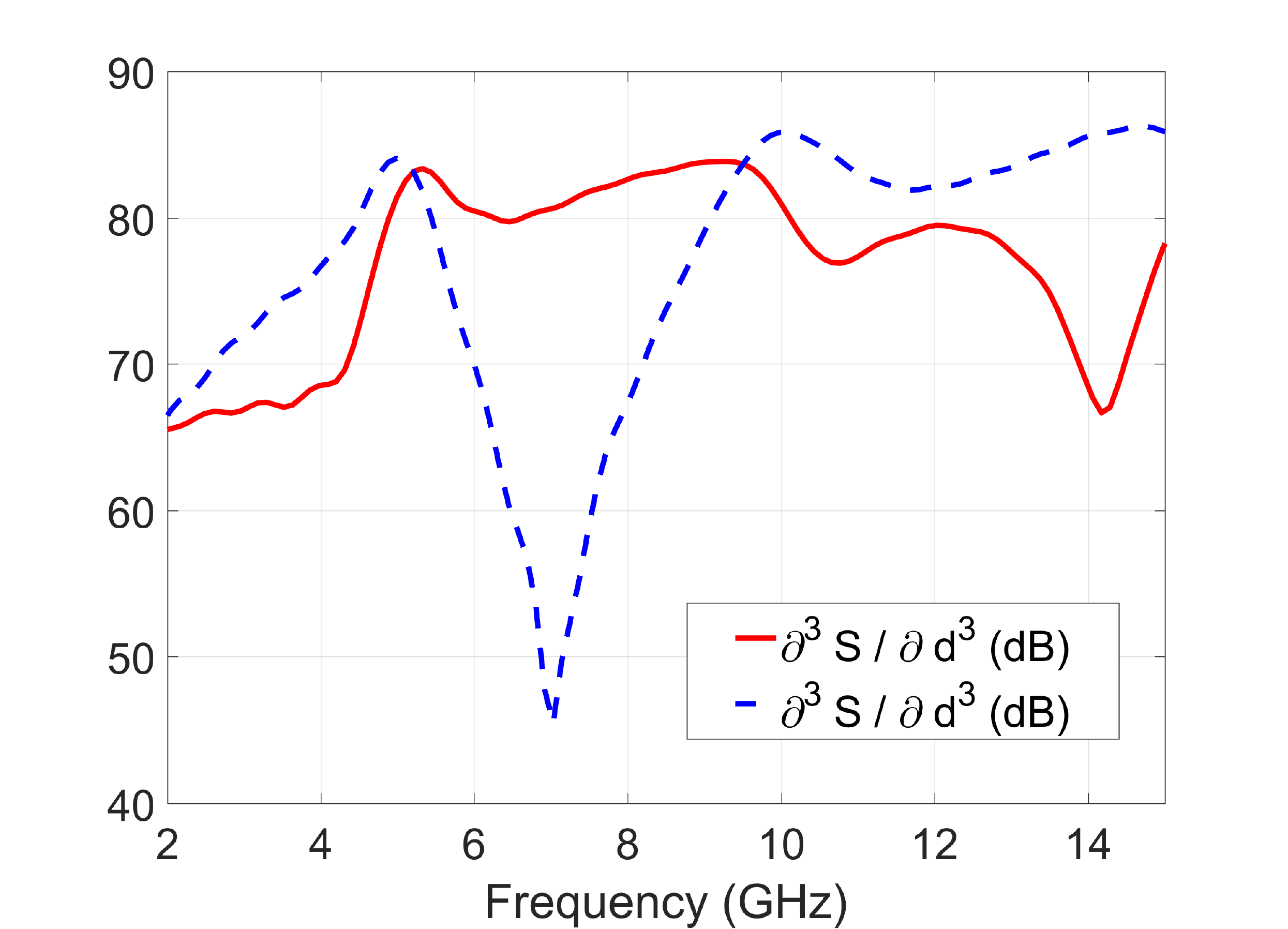}}
		\caption{Numerical sensitivities $\partial^3 S_{11}, S_{21} / \partial w_1\partial w_2 \partial \epsilon_r$ (a) and $\partial^3 S_{11}, S_{21} / \partial d^3$ (b) of the microstrip low-pass filter, computed via FDTD with the multi-complex step derivative (MCSD) approximation.}
		\label{filter3rd}
\end{figure}

\textcolor{black}{
Finally, the application of these high-order field derivatives to the parametric modelling of output functions of interest is presented. For example, a parametric model of $S_{21}$ with respect to substrate thickness $d$, can be derived from a Taylor expansion around the nominal value $d_0 = 0.794$ mm:
\begin{equation}
S_{21} (d_0) = \sum_{n=0}^\infty\frac{S_{21}^{(n)} (d_0)}{n!} (x-d_0)^n \\ 
\label{Taylor expansion}
\end{equation}
}
\begin{figure}[]
	 \centering
		\subfigure[]{\includegraphics[width=9cm]{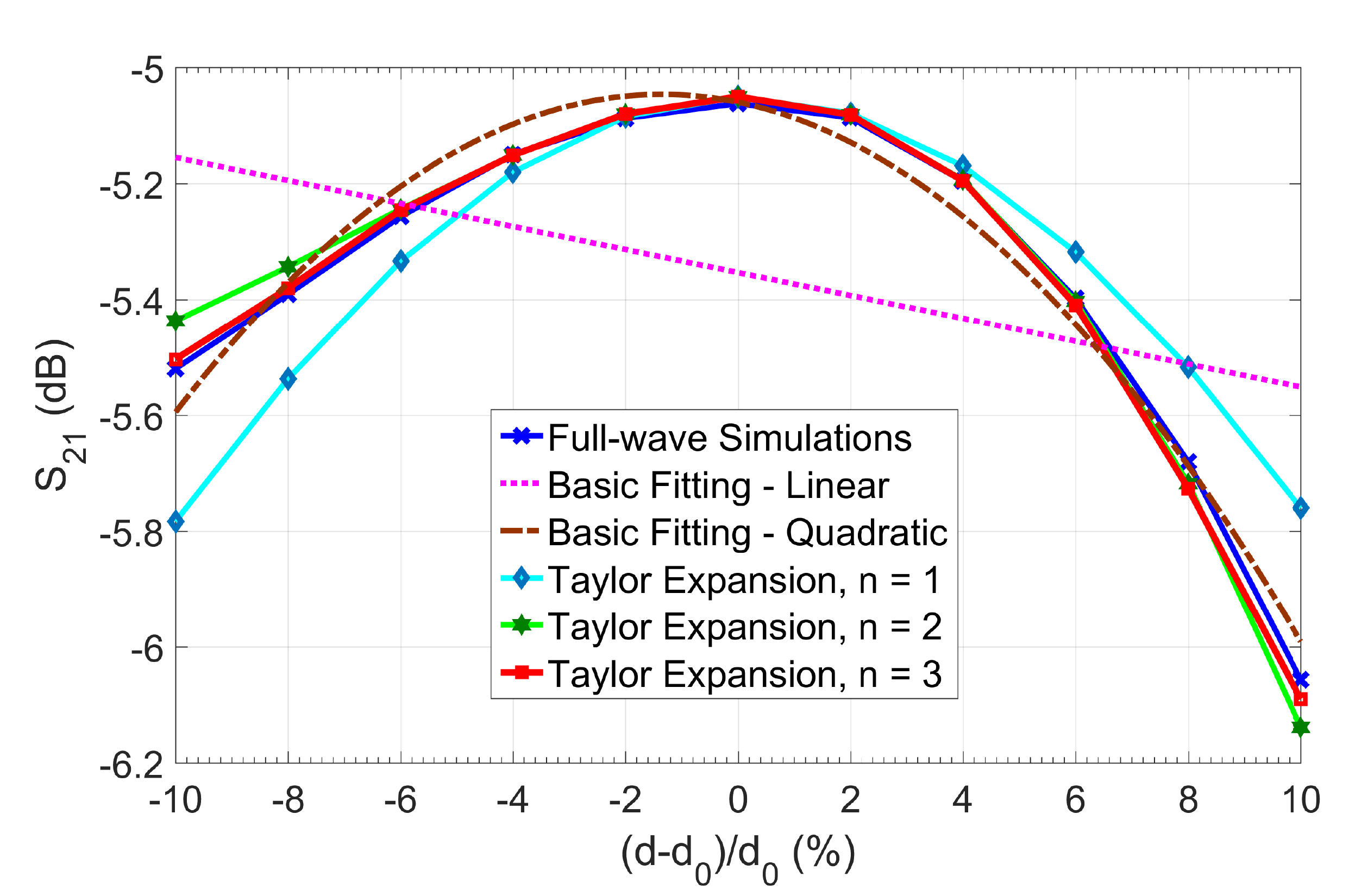}}
		\subfigure[]{\includegraphics[width=9cm]{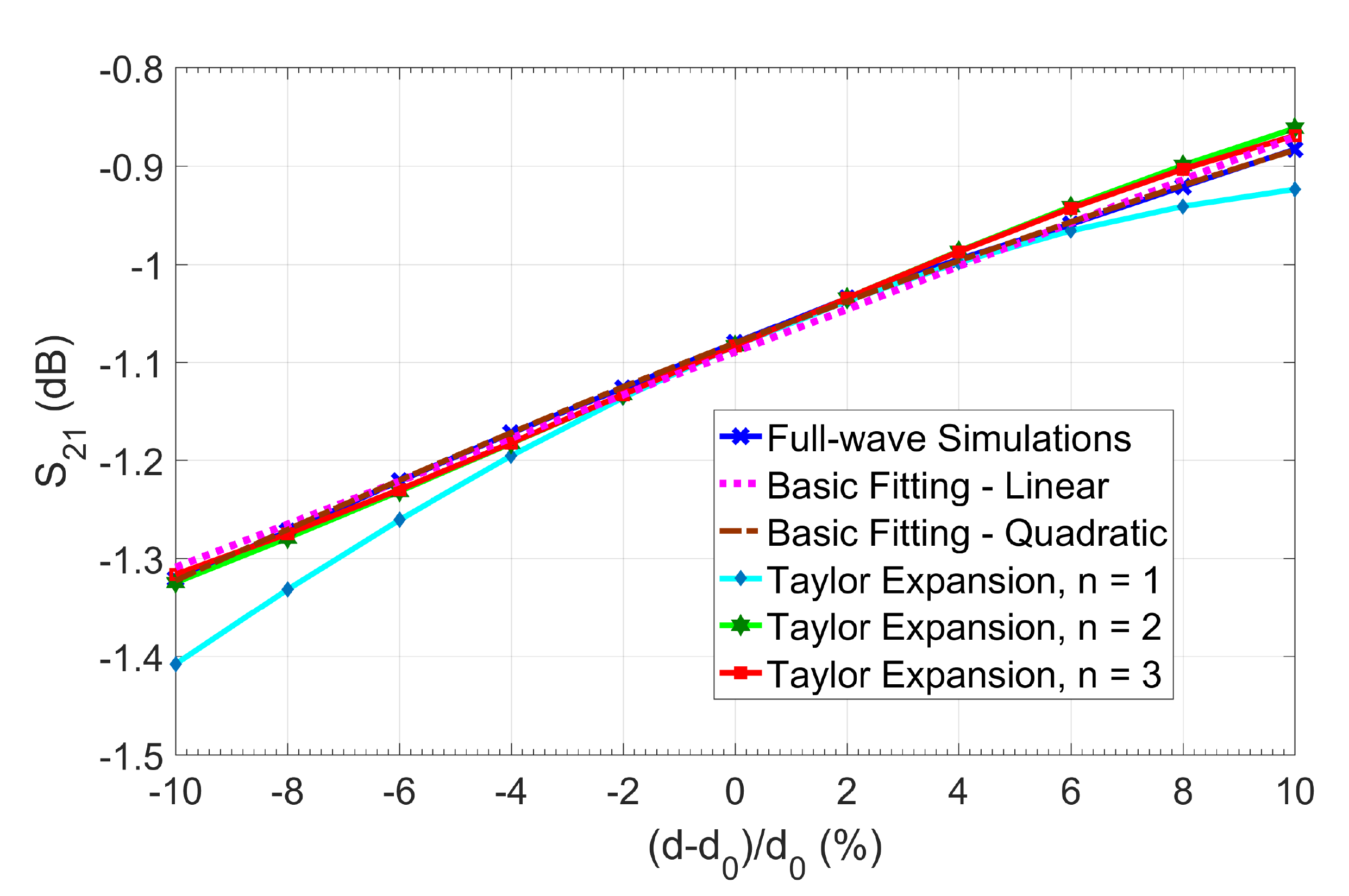}}
		\caption{\textcolor{black}{The insertion loss $S_{21}$ of the microstrip filter with substrate thickness varying from $-10$ to $10\%$ of the nominal value $d_0 = 0.753$ mm at frequency of (a) 4.2 GHz and (b) 2.2 GHz. The effective range of parametric model is increased as higher order Taylor expansion terms are included. At frequency of 4.2 GHz, linear and quadratic line fitting of full-wave simulations data points fail to model the non-linear relation between $S_{21}$and $d$. }}
\label{Taylor}
\end{figure}
\textcolor{black}{
Fig.\ref{Taylor} shows the $S_{21}$ estimated by a Taylor expansion with thickness $d$ varying from $0.9d$ to $1.1d$, at frequencies of 4.2 GHz and 2.2 GHz. Full-wave simulation results are shown along with the parametric model. It is found that the insertion loss is more sensitive to thickness at 4.2 GHz, which is at the edge of the first pass-band of the filer, and higher order derivatives become important in building an accurate parametric model. 
}
\textcolor{black}{
Notably, the effective bandwidth of the model is increased from $4\%$ to $18\%$ at 4.2 GHz by increasing the order of the Taylor expansion. This shows the importance of high-order derivative computations that is enabled by MCSD and the iterative CSD scheme. The wide effective range allows users to estimate the performance of electromagnetic design with varying parameters.
}
\section{Conclusions}
\begin{table}[]
\centering
\caption{Comparison on FDTD-based Field Derivative Computation Methods }
\label{mcsd-cfd-avmcompare}
\begin{tabular}{|p{1.8cm}|p{1.8cm}|p{1.8cm}|p{1.8cm}|}
\hline
 & MCSD-FDTD & CFD-FDTD & AVM-FDTD\\
 \hline
Jacobian Complexity & $O(N)$ & $O(N)$ & $O(1)$ \\
\hline
Hessian Complexity& $O(N^2)$ & $O(N^2)$ & $O(N)$ \\
\hline
Accuracy& Machine accuracy&Subtraction-error, step-size needs testing & Comparable to CFD\\
\hline
Implementation & Simple once complex numbers are defined & Simple; runs standard FDTD & Non-trivial \\
\hline
Advantage & Simple implementation, scalability, accuracy&  Simple implementation  & Efficient for multiple parameters 
\\
\hline
\end{tabular}
\end{table}
This paper presented a general framework for accurate computation of high-order field derivatives with respect to design parameters
in FDTD simulations. Indeed, the accuracy and robustness of the MCSD approximation and the relative simplicity of its implementation in FDTD
(independently of the order of derivative needed), just by including an appropriately defined class of multi-complex numbers, are the most significant advantages of this method. In a nutshell, MCSD-FDTD offers the versatility of finite-differences, as it allows for the computation of field derivatives of any order with FDTD, while its accuracy significantly surpasses that of finite-difference methods, especially for small perturbation steps.  

\textcolor{black}{A detailed analysis of the computational cost of CSD and MCSD was presented.
Based on this analysis, first-order derivatives with respect to multiple design parameters can be more efficiently computed by applying 
CSD-FDTD to  each parameter rather than a full MCSD-FDTD analysis. For high-order derivatives, CSD-FDTD can be iteratively applied in a marching-in-order scheme to alleviate the computational overhead of MCSD-FDTD. This is being said, MCSD-FDTD does provide a complete framework to compute field derivatives of any order. To summarize, a comparison between the proposed method and two popular alternatives, CFD and AVM respectively, is presented in Table \ref{mcsd-cfd-avmcompare}. Notably, AVM has so far been formulated
towards the computation of  derivatives of a given, field-based  objective function rather than  field derivatives themselves.}


While this paper focused on  FDTD, the underlying theory of complex step derivative approximations can clearly be combined with other full-wave and equivalent circuit analysis techniques (in frequency and time-domain), employed in electromagnetic design. 

\bibliographystyle{IEEEtran}
{\footnotesize
\bibliography{IEEEfull,../../mybibfile}
}
\appendices
\end{document}